\documentclass[a4paper, 10pt]{amsart}
\usepackage{times,latexsym,amssymb}
\usepackage{amsmath}
\usepackage{amsthm}

\allowdisplaybreaks
\newtheorem{theorem}{Theorem}
\newtheorem*{theorem-plain}{Theorem}
\newtheorem{lemma}[theorem]{Lemma}
\newtheorem{definition}[theorem]{Definition}
\newtheorem{proposition}[theorem]{Proposition}

\newtheorem{remark}[theorem]{Remark}

\numberwithin{theorem}{section}
\numberwithin{equation}{section}

\newcommand{\mint}{- \mskip-19,5mu \int}

\def\N{\mathbb{N}}
\def\R{\mathbb{R}}

\def\dz{\,dz}
\def\dx{\,dx}
\def\dy{\,dy}
\def\dt{\,dt}
\def\ds{\,ds}

\def\G{{\boldsymbol G}}

\def\V{{\boldsymbol V}}

\def\<{\langle}
\def\>{\rangle}

\def\nn{\nonumber}

\DeclareMathOperator{\osc}{osc}
\DeclareMathOperator{\Div}{div}

\DeclareMathOperator\VMO{VMO}
\renewcommand{\epsilon}{\varepsilon}
\renewcommand{\rho}{\varrho}

\begin{document}
\title[Global Calder\'on \& Zygmund theory]{Global Calder\'on \& Zygmund theory for nonlinear parabolic systems}
\date{\today}
\author[V. B\"ogelein]{Verena B\"{o}gelein}
\address{Verena B\"ogelein\\Department Mathematik, Universit\"at
Erlangen--N\"urnberg\\ Cauerstra\ss e 11, 91058 Erlangen, Germany}
\email{boegelein@math.fau.de}


\begin{abstract} 
We establish a global Calder\'on \& Zygmund theory for solutions of a huge class of nonlinear parabolic systems whose model is the inhomogeneous parabolic $p$-Laplacian system 
\begin{equation*}
\left\{
\begin{array}{cc}
	\partial_t u - \Div (|Du|^{p-2}Du) = \Div (|F|^{p-2}F)
	&\mbox{in $\Omega_T:=\Omega\times(0,T)$} \\[5pt]
    u=g
    &\mbox{on $\partial\Omega\times(0,T)\cup \overline\Omega\times\{0\}$}
\end{array}
\right.
\end{equation*}
with given functions $F$ and $g$. Our main result states that the spatial gradient of the solution is as integrable as the data $F$ and $g$ up to the lateral boundary of $\Omega_T$, i.e.
\begin{equation*}
	F,Dg\in L^q(\Omega_T),\ \
	\partial_t g\in L^{\frac{q(n+2)}{p(n+2)-n}}(\Omega_T)
	\quad\Rightarrow \quad
	Du\in L^q(\Omega\times(\delta,T))
\end{equation*}
for any $q>p$ and $\delta\in(0,T)$, together with quantitative estimates.
This result is proved in a much more general setting, i.e.  for  asymptotically regular parabolic systems. 
\end{abstract}

\maketitle
%*****************************************************************
\tableofcontents
%\newpage
%*****************************************************************
\section{Introduction}

We establish a global Calder\'on \& Zygmund theory, i.e. $L^q$ regularity results for the (spatial) gradient of the solution, for a huge class of nonlinear parabolic systems.
As model problem we consider the inhomogeneous parabolic $p$-Laplacian system 
\begin{equation}\label{para-plap-inhom}
	\partial_t u - \Div (|Du|^{p-2}Du) = \Div (|F|^{p-2}F)
\end{equation}
with a given function $F$. The scope of the theory is to ensure that the solution $u$ is as good integrable as the inhomogeneity $F$, i.e. the local version reads as
\begin{equation*}
	F\in L^q_{\rm loc} \Rightarrow Du\in L^q_{\rm loc}
	\quad \mbox{for any $q>p$.}
\end{equation*}
For the stationary, elliptic case such a result has been obtained by 
Iwaniec \cite{Iwaniec:1983} in the scalar case ($N=1$) and by 
DiBenedetto \& Manfredi \cite{DiBenedetto-Manfredi:1993} in the vectorial case ($N>1$).
This result has been extended to elliptic equations with VMO coefficients by Kinnunen \& Zhou \cite{Kinnunen-Zhou-1999}. The idea how to treat more general non-linear elliptic equations is due to 
Caffarelli \& Peral \cite{Caffarelli-Peral:1998}. Based on this technique a general 
Calder\'on \& Zygmund theory for elliptic equations and systems including the case of non-standard growth conditions has been obtained by Acerbi \& Mingione \cite{Acerbi-Mingione:2005}. These results were also achieved in a global version for homogeneous Dirichlet boundary data, i.e. $u=0$ on the boundary by Kinnunen \& Zhou \cite{Kinnunen-Zhou-2001} for $C^{1,\alpha}$ domains and by Byun \& Wang \cite{Byun-Wang:2008-2} in Reifenberg domains. A global Calder\'on \& Zygmund theory for minimizers of integral functionals with non-homogeneous boundary data has been established by Kristensen \& Mingione
\cite{Kristensen-Mingione:2010}.

First of all it was not clear how to prove a similar result in the parabolic setting. Thereby, the main obstruction was the non-homogeneous scaling behavior of the problem with respect to space and time in the case $p\not=2$, in the sense that the solution multiplied by a constant is in general not anymore a solution. As a consequence all basic estimates -- such as reverse H\"older inequalities -- become inhomogeneous and the use of the maximal function becomes delicate.
A first result in this direction for the special case when $F\in L^\infty$ has been obtained by Misawa \cite{Misawa:2006}. 
The local Calder\'on \& Zygmund theory for the time dependent, parabolic case together with a quantitative estimate has been achieved by Acerbi \& Mingione \cite{Acerbi-Mingione:2007} and for more general non-linear systems by Duzaar \& Mingione \& Steffen \cite{Duzaar-Mingione-Steffen:2011}. 
In the setting of obstacle problems we refer to \cite{Boegelein-Duzaar-Mingione:2011, Scheven:2011}. 
The main idea to overcome the difficulties in the parabolic case was to use DiBenedettos intrinsic geometry together with a maximal function free approach.
Up to now a global result was available only for the case $p=2$ and homogeneous Dirichlet boundary data, i.e. $u=0$ on the parabolic boundary, see Byun \& Wang \cite{Byun-Wang:2008}. 
As already mentioned before, the case $p\not=2$ is much more involved since the scaling of the underlying system is non-homogeneous.
Therefore, our main purpose in this paper is to establish the following global Calder\'on \& Zygmund theorem for parabolic systems of $p$-Laplacian type.

\begin{theorem-plain}
Let $p>\frac{2n}{n+2}$ and suppose that
$$u\in C^0\big([0,T];L^{2} (\Omega,\R^N)\big)\cap L^p\big(0,T;W^{1,p}(\Omega,\R^N)\big)$$ 
is a weak solution to the Cauchy-Dirichlet problem 
\begin{equation*}
\left\{
\begin{array}{cc}
    \partial_t u - \Div \big(c(z)|Du|^{p-2}Du\big) = \Div(|F|^{p-2}F) + f
    &\mbox{in $\Omega_T$} \\[5pt]
    u=g
    &\mbox{on $\partial\Omega\times(0,T)\cup \overline\Omega\times\{0\}$}
\end{array}
\right.
\end{equation*}
where $\Omega_T:=\Omega\times(0,T)$, $\Omega$ is a bounded $C^1$ domain and $c\colon\Omega_T\to [\nu,L]$, $0<\nu\le L$ is VMO with respect to $x$ (see \eqref{VMO-c} below) and measurable with respect to $t$.
Further, assume that 
\begin{equation*}
	F,Dg\in L^q(\Omega_T,\R^{Nn}),\qquad
	f,\partial_t g\in L^{\frac{q(n+2)}{p(n+2)-n}}(\Omega_T,\R^N)
\end{equation*}
for some $q>p$.
Then
\begin{equation*}
	D u\in L^{q}\big(\Omega\times(\delta,T),\R^{Nn}\big)
	\quad\mbox{for any $\delta\in(0,T)$}
\end{equation*}
and there holds the quantitative $L^q$-estimate \eqref{CZ-est} below.
\end{theorem-plain}

We note that for the sake of brevity we only consider the lateral boundary regularity, keeping in mind that one could extend the theory to the initial boundary.

Furthermore, we will obtain our result not only for the pure $p$-Laplacian system, but also for a much larger class of parabolic systems, the so called asymptotically regular systems. By this we mean parabolic systems of the type 
\begin{equation*}
    \partial_t u - \Div a(z,u,Du) = \Div(|F|^{p-2}F) + f,
\end{equation*}
where the vector field $a\colon\Omega_T\times\R^N\times\R^{Nn}\to\R^{Nn}$ is asymptotically regular in the sense that 
\begin{equation}\label{asy-intro}
    \lim_{|\xi|\to\infty}
    \frac{a(z,u,\xi) - b(z,u,\xi)}{|\xi|^{p-1}}
    =
    0
\end{equation}
holds uniformly with respect to $z\in\Omega_T$ and $u\in\R^N$ for some regular vector field $b\colon\Omega_T\times\R^N\times\R^{Nn}\to\R^{Nn}$. The notion regular will be specified later. We will prove that solutions to the associated parabolic Cauchy-Dirichlet problem with vector field $a$ are almost as integrable as solutions to the same Cauchy-Dirichlet problem with vector field $b$. A Calder\'on \& Zygmund theory for this kind of systems is new even in the interior situation. 
The crucial point here is that \eqref{asy-intro} is essentially the only assumption imposed on the vector field $a$ itself. In particular this means that there is {\it no assumption for small $\xi$}.
%For example, if $b$ has linear growth, is independent of $u$, satisfies a VMO-condition with respect to $x$ and is measurable with respect to $t$, then the result of the last theorem holds whenever $q<2+\frac4n+\epsilon$ with some $\epsilon>0$ depending only on the structural data. The upper bound on $q$ in this case is natural, since even for the associated homogeneous problem it is not true that $Du\in L^q$ for any $q\ge p$.
For the precise statements and applications of the general result we refer to Section \ref{sec:results-CZ}.

The notion of asymptotic regular problems was introduced in the elliptic framework by Chipot \& Evans \cite{Chipot-Evans:1986}. They proved Lipschitz regularity of minimizers to integral functionals $F(v)=\int_\Omega f(Dv)\dx$ with an integrand satisfying $D^2f(\xi)\to A$ when $|\xi|\to\infty$ for some elliptic bilinear form $A$ on $\R^{Nn}$. More general integrands were treated later by 
Giaquinta \& Modica \cite{Giaquinta-Modica:1986} and Raymond \cite{Raymond:1991} and the case of higher order functionals has been considered by Schemm \cite{Schemm:2009}.
Finally, for a result in the context of Orlicz spaces we refer to Diening \& Stroffolini \& Verde \cite{Diening-Stroffolini-Verde:2011}. Global Morrey and Lipschitz regularity results have been obtained by 
Foss \cite{Foss:2008} and Foss \& Passarelli di Napoli \& Verde \cite{Foss-Passarelli-Verde:2008}.
A local Calder\'on \& Zygmund theory and partial Lipschitz regularity for asymptotically regular elliptic systems and minimizers has been developed by Scheven \& Schmidt \cite{Scheven-Schmidt:2009, Scheven-Schmidt:2010}.
These results provide a huge class of integral functionals with (partial) Lipschitz minimizers, which is much larger than the well known class of quasi-diagonal structure.
In the evolutionary framework, the global Lipschitz regularity up to the lateral boundary of solutions to parabolic asymptotically $p$-Laplacian systems has been proved in \cite{Boegelein:2012-Lipschitz, Boegelein-habil}.

\section{Statement of the results}\label{sec:results-CZ}

We fix $n\ge 2$, $N\ge 1$ and a growth exponent $p>\frac{2n}{n+2}$.
With $\Omega_T:=\Omega\times(0,T)\subset\R^{n+1}$ denoting the space-time cylinder over a bounded domain $\Omega\subset\R^n$ and $T>0$ we suppose that 
$a\colon\Omega_T\times\R^N\times\R^{Nn}\to \R^{Nn}$ is a Carath\'eodory vector field and $f,g\colon\Omega_T\to\R^N$ and $F\colon\Omega_T\to\R^{Nn}$ are some given functions.
Then, we are interested in parabolic Cauchy-Dirichlet problems of the type 
\begin{equation}\label{Cauchy-Dirichlet}
\left\{
\begin{array}{cc}
    \partial_t u - \Div a(z,u,Du) = \Div(|F|^{p-2}F) + f
    &\qquad\mbox{in $\Omega_T$} \\[5pt]
    u=g
    &\qquad\mbox{on $\partial_{\mathcal P}\Omega_T$.}
\end{array}
\right.
\end{equation}
Note that the solution $u\colon\Omega_T\to\R^N$ is possibly vector valued and that by $\partial_{\mathcal P}\Omega_T:=\partial\Omega\times(0,T)\cup \overline\Omega\times\{0\}$ we denote the parabolic boundary of $\Omega_T$. The only assumption on the vector field $a$ will be that it is uniformly asymptotically related to a more regular Carath\'eodory vector field $b\colon\Omega_T\times\R^N\times\R^{Nn}\to \R^{Nn}$, in the sense that 
\begin{equation}\label{def_asymp-CZ}
    \lim_{|\xi|\to\infty}
    \frac{a(z,u,\xi) - b(z,u,\xi)}{|\xi|^{p-1}}
    =
    0
%    \qquad\mbox{with $p>\frac{2n}{n+2}$}
\end{equation}
holds uniformly with respect to $z\in\Omega_T$ and $u\in\R^N$. The admissible classes of vector fields $b$ will be specified in Sections \ref{sec:asym-plap} -- \ref{sec:general-setting} below.
Concerning the regularity of the boundary data, i.e. of $\partial\Omega$ and $g$, we shall assume that 
%$g\colon \overline{\Omega}_T \to \R^N$ is a continuous function such that
\begin{equation}\label{data-CZ}
	\partial\Omega \text{ is } C^{1},\quad
	Dg\in L^q(\Omega_T,\R^{Nn}),\quad
	\partial_t g\in L^{\frac{q(n+2)}{p(n+2)-n}}(\Omega_T,\R^N)
\end{equation}
for some $q>p$. Moreover, the functions $F$ and $f$ on the right-hand side of \eqref{Cauchy-Dirichlet} are supposed to satisfy
\begin{equation}\label{inhom}
	F\in L^q(\Omega_T,\R^{Nn}),\quad
	f\in L^{\frac{q(n+2)}{p(n+2)-n}}(\Omega_T,\R^N).
\end{equation}
Next, we specify the notion of a weak solution
to the Cauchy-Dirichlet problem \eqref{Cauchy-Dirichlet}.
\begin{definition}\label{def:weak-solution-CZ}
A map 
$$
	u \in C^0\big([0,T];L^2(\Omega,\R^{N})\big)\cap L^p\big(0,T;W^{1,p}(\Omega,\R^{N})\big)
$$ 
is called a (weak) solution to the parabolic Cauchy-Dirichlet problem \eqref{Cauchy-Dirichlet}
if and only if
\begin{align}\label{weak}
	\int_{\Omega_T}
	u\cdot\varphi_t - \langle a(z,u,Du), D\varphi\rangle \dz
	=
	\int_{\Omega_T}
	\langle |F|^{p-2}F, D\varphi\rangle - f\cdot\varphi \dz
\end{align}
whenever $\varphi\in C_0^\infty(\Omega_T,\R^N)$, and the following boundary conditions are satisfied:
\begin{equation*}
	u(\cdot,t) - g(\cdot,t) \in W^{1,p}_0(\Omega,\R^N)
	\qquad\text{for a.e. } t\in(0,T)
\end{equation*}
and
$$
  \lim_{h \downarrow0 }\frac{1}{h} \int_{0}^{h} \int_\Omega |u(x,t) - g(x,0)|^2 \dx\dt =0
  \,.
$$
\end{definition}

In the following we will provide a local and global Calder\'on \& Zygmund theory for a huge class of non-linear parabolic systems. It applies to (asymptotically) parabolic $p$-Laplacian systems as well as to their counterparts without a quasi-diagonal structure.
With this respect it will be more convenient to work with the following equivalent notion of asymptotic regularity.
\begin{remark}\label{rem:asymp}\upshape
The vector field $a$ is asymptotically regular in the sense of \eqref{def_asymp-CZ} if and only if 
\begin{equation}\label{def_asymp-CZ-}
    |a(z,u,\xi) - b(z,u,\xi)|
    \le
    \omega(|\xi|)(1+|\xi|)^{p-1},
    \quad\forall\, z\in\Omega_T, u\in\R^N, \xi\in \R^{Nn}
\end{equation}
holds for some bounded function $\omega\colon[0,\infty)\to [0,\infty)$ with
\begin{equation}\label{omega-CZ}
    \lim_{s\to \infty} \omega(s)
    = 0.
\end{equation}
\vspace{-0.7cm}

\hfill$\Box$
\end{remark}

\subsection{Asymptotically $p$-Laplacian systems}\label{sec:asym-plap}
The first result concerns a global Calder\'on \& Zygmund theory for the parabolic $p$-Laplacian system, or more generally for asymptotically parabolic $p$-Laplacian systems. More precisely, in \eqref{def_asymp-CZ} we assume that the vector field $b$ is of the form
\begin{equation}\label{b-1}
    b(z,u,\xi)
    \equiv
    b(z,\xi)
    =
    c(z) |\xi|^{p-2}\xi
    \qquad\mbox{for $z\in \Omega_T, \xi\in\R^{Nn}$,}
\end{equation}
where $c\colon\Omega_T\to [\nu,L]$, $0<\nu\le L$ is measurable with respect to $t$ and VMO with respect to $x$, i.e. 
\begin{equation}\label{VMO-c}
   	\V(\rho)
	:=
	\sup_{t\in(0,T)}
	\sup_{x_o\in\Omega}
    	\sup_{0<r\le\rho}\mint_{B_r(x_o)\cap \Omega} \big|c(x,t)-\big(c(\cdot,t)\big)_{B_r(x_o)\cap \Omega} \big|\dx
    	\to 0
\end{equation}
as $\rho\downarrow 0$. Then, we have the following

\begin{theorem}\label{thm:main-CZ1}
Let $q>p$ and suppose that
$$u\in C^0\big([0,T];L^{2} (\Omega,\R^N)\big)\cap L^p\big(0,T;W^{1,p}(\Omega,\R^N)\big)$$ 
is a weak solution to the parabolic Cauchy-Dirichlet problem \eqref{Cauchy-Dirichlet}, under the assumptions \eqref{def_asymp-CZ} -- \eqref{inhom} with \eqref{b-1} and \eqref{VMO-c}.
Then
\begin{equation*}
	D u\in L^{q}\big(\Omega\times(\delta,T),\R^{Nn}\big)
	\quad\mbox{for any $\delta\in(0,T)$.}
\end{equation*}
Moreover, there exists $R_o>0$ such that for any $z_o\in \overline\Omega\times(0,T)$ and any parabolic cylinder $Q_{2R}(z_o)\subset\R^n\times(0,T)$ with $R\in(0,R_o]$ there holds
\begin{align}\label{CZ-est}
	& \mint_{Q_{R}(z_o)\cap \Omega_T} |D u|^{q} \dz \\
	&\phantom{mm}\le
	c \Bigg[\bigg(\mint_{Q_{2R}(z_o)\cap \Omega_T}|D u|^p\dz \bigg)^{\frac{q}{p}} +
	\mint_{Q_{2R}(z_o)\cap \Omega_T}
	(|Dg| + |F|)^q \dz \nn\\
	&\phantom{mmmm}+
	R^{-(n+2)} \bigg(\int_{Q_{2R}(z_o)\cap \Omega_T}
	(|g_t| + |f| )^{\frac{q(n+2)}{p(n+2)-n}} \dz 
	\bigg)^{1+\frac{p}{np+p-n}} +1
	\Bigg]^{d_{CZ}}, \nn
\end{align}
where we have abbreviated
\begin{equation}\label{def-d}
	d_{CZ}
	:=
	d-\tfrac{p}{q}(d-1)
	\quad\mbox{with }\
	d:=
	\begin{cases}
	\frac{p}{2} & \text{if } p\ge 2,\\[5pt]
	\frac{2p}{p(n+2) - 2n} & \text{if } p<2.
	\end{cases}
\end{equation}
Note that the constant $c$ depends on $n,N,p,\nu,L,q,\omega(\cdot), \V(\cdot),\partial\Omega$ and 
$R_o$ depends on 
$n,$ $N,p,\nu,L,q, \|\omega\|_\infty, \V(\cdot),\partial\Omega$, where $\omega(\cdot)$ is from Remark \ref{rem:asymp}.
\end{theorem}

As we already mentioned in the introduction, the global result is new even for the model case of the parabolic $p$-Laplacian, i.e. $a(z,u,\xi)=|\xi|^{p-2}\xi$.
Since the result for asymptotically parabolic $p$-Laplacian systems is new also for the interior situation, we state the simpler quantitative estimate for this case.

\begin{remark}\label{rem:interior}\upshape
If $Q_{2R}(z_o)$ is an interior cylinder, i.e. if
$Q_{2R}(z_o)\subset\Omega_T$, the quantitative $L^q$-estimate from Theorem \ref{thm:main-CZ1} holds without the boundary terms $|Dg|$ and $|\partial_t g|$. More precisely, in this case we have
\begin{align}\label{CZ-est-int}
	\mint_{Q_{R}(z_o)} |D u|^{q}\dz
	&\le
	c\, \Bigg[\bigg(\mint_{Q_{2R}(z_o)}|D u|^p\dz \bigg)^{\frac{q}{p}} +
	\mint_{Q_{2R}(z_o)} |F|^q \dz \\
	&\phantom{mm}+
	R^{-(n+2)}
	\bigg(\int_{Q_{2R}(z_o)} |f|^{\frac{q(n+2)}{p(n+2)-n}} \dz 
	\bigg)^{1+\frac{p}{np+p-n}} + 1
	\Bigg]^{d_{CZ}} \nn
\end{align}
provided $R\in(0,R_o]$. 
Note that now the constant $c$ depends on $n,N,p,\nu,L,q,\omega(\cdot), \V(\cdot)$ and 
$R_o$ depends on 
$n,N,p,\nu,L,q, \|\omega\|_\infty, \V(\cdot)$.
This can easily be seen from the proof, since one can omit the transformation to the model situation on $Q_{2R}^+$.
\hfill$\Box$
\end{remark}

\begin{remark}\upshape
Here, we comment on the exponents and constants appearing in Theorem~\ref{thm:main-CZ1}. First of all, $d$ is called the parabolic scaling deficit which naturally appears in estimates for parabolic problems with $p$-growth in the diffusion term, since the space-time scaling is non-homogeneous unless $p=2$. Indeed, in the case $p=2$ it happens that $d=1$ and otherwise it holds that $d>1$. 
The exponent $d_{CZ}$ is the improved scaling deficit appearing in the Calder\'on \& Zygmund estimate \eqref{CZ-est}, respectively \eqref{CZ-est-int}. If $p=2$ we have $d_{CZ}=1$ and $1<d_{CZ}<d$ otherwise.
Moreover, we have
$$
	d_{CZ} \downarrow 1
	\quad\mbox{when $q\downarrow p$}
$$
which seems to be natural. Furthermore, in contrast to \cite[Remark 3]{Acerbi-Mingione:2007}, the constants in Theorem \ref{thm:main-CZ1} and Remark \ref{rem:interior} remain stable when $q\downarrow p$ without any stabilizing procedure.

The integrability exponent $\frac{q(n+2)}{p(n+2)-n}$ of the non-divergence form inhomogeneity $f$ and of $\partial_t g$ arises from an application of the parabolic Sobolev-Poincar\'e inequality, see Section \ref{sec:comparison-CZ}. A straightforward application of H\"older's inequality would yield the strictly larger exponent $\frac{q}{p-1}$ instead.
\hfill$\Box$
\end{remark}

\subsection{General non-linear parabolic systems}\label{sec:asym-gen}
Here we consider more general vector fields $b\colon\Omega_T\times\R^N\times\R^{Nn}\to\R^{Nn}$ which are not necessarily of quasi-diagonal structure. More precisely, we assume that $b$ and $\partial_\xi b$ are Carath\'eodory maps satisfying 
\begin{align}\label{assump-b}
\left\{
\begin{array}{c}
	|b(z,u,\xi)| + \big(\mu^2+|\xi|^2\big)^{\frac12}|\partial_\xi b(z,u,\xi)|
	\le
	L\,\big(\mu^2+|\xi|^2\big)^{\frac{p-1}{2}}\\[5pt]
	\langle\partial_\xi b(z,u,\xi)\xi_o,\xi_o\rangle
	\ge
	\nu\,\big(\mu^2+|\xi|^2\big)^{\frac{p-2}{2}} |\xi_o|^2\\[5pt]
	|b(z,u,\xi) - b(z,u_o,\xi)|
	\le
	\theta\big(|u-u_o|\big)\big(\mu^2+|\xi|^2\big)^{\frac{p-1}{2}}
\end{array}
\right.
\end{align}
whenever $z\in\Omega_T$, $u,u_o\in\R^N$, $\xi,\xi_o\in\R^{Nn}$, for some parameters $0<\nu\le 1\le L$ and $\mu\in[0,1]$. In the case $p<2$ and $\mu=0$ we certainly assume that $\xi\not=0$ in \eqref{assump-b}$_2$. Thereby, $\theta\colon[0,\infty)\to [0,1]$ is a nondecreasing, concave modulus of continuity with $\lim_{s\downarrow 0}\theta(s)=0=\theta(0)$.
With respect to the dependence on $t$ the map $b$ is only assumed to be measurable, while with respect to $x$ we impose a VMO-condition; more precisely we assume that
$x\mapsto b(x,t,u,\xi)/(1+|\xi|)^{p-1}$ fulfills the following $\VMO$-condition uniformly in $t,u$ and $\xi$:
\begin{equation}
  \label{VMO-b}
    \big|b(x,t,u,\xi)- \big( b(\cdot,t , u,\xi)\big)_{\Omega\cap B_r (x_o)}\big|
    \le
    {\bf v}_{x_o}(x,r)\big(\mu^2+|\xi|^2\big)^{\frac{p-1}{2}}
    \quad\forall\, x\in B_r(x_o)
\end{equation}
whenever $x_o\in\Omega$, $t\in(0,T)$,
$u\in\R^N$, $\xi\in \R^{Nn}$ and $r\in(0,\rho_o]$,  where $\rho_o>0$ and ${\bf v}_{x_o} \colon \R^n\times [0,\rho_o]\to
[0,2L]$ are  bounded functions satisfying
\begin{equation}
  \label{VMO-b-}
    \V(\rho ):=\sup_{x_o\in\Omega}\,\sup_{0<r\le\rho}
    \mint_{\Omega\cap B_r (x_o)}{\bf v}_{x_o}(x,r)\, dx
    \to 0
    \quad\mbox{as } \rho\downarrow 0.
\end{equation}
Here, we used the short-hand notation
$$
  \big(b(\cdot ,t,u,\xi)\big)_{\Omega\cap B_r (x_o)}:=
  \mint_{\Omega\cap B_r (x_o)} b(y,t,u,\xi)\, dy
$$
for the partial means -- i.e. the means with respect to $x$ -- of
the vector-field $b$ for fixed arguments $(t,u,\xi)\in(0,T)\times\R^N\times\R^{Nn}$.
Then, in the interior situation we get the following result.

\begin{theorem}\label{thm:main-CZ3}
There exists $\epsilon=\epsilon(n,N,p,\nu,L)>0$ such that the following holds: whenever $q\in(p,p+\frac4n+\epsilon]$ and 
$$u\in C^0\big([0,T];L^{2} (\Omega,\R^N)\big)\cap L^p\big(0,T;W^{1,p}(\Omega,\R^N)\big)$$ 
is a weak solution to the parabolic Cauchy-Dirichlet problem \eqref{Cauchy-Dirichlet}, under the assumptions \eqref{def_asymp-CZ} -- \eqref{inhom} with \eqref{assump-b} -- \eqref{VMO-b-} and 
\begin{equation*}
	b(z,u,\xi)\equiv b(z,\xi),
\end{equation*}
then we have
\begin{equation*}
	D u\in L^{q}_{\rm loc}\big(\Omega_T,\R^{Nn}\big).
\end{equation*}
Moreover, the quantitative estimate \eqref{CZ-est-int} holds.
\end{theorem}

It is not clear if the previous result can be extended up to the lateral boundary. The only obstruction thereby, is to obtain up-to-the-boundary a priori estimates. More precisely, the difficulty relies in the fact that at the lateral boundary two directions -- the normal spatial and the time direction -- have to be recovered from the parabolic system. 
On the other hand, in the particular case $p=2$ it is possible to prove up-to-the-boundary $L^2$-estimates for the first time-derivative by the use of second finite differences which in turn allow for boundary a priori estimates; see Lemma \ref{lem:chi-bound-2} and \cite[Chapter 4]{Boegelein-Duzaar-Mingione:2010-boundary-II} for the complete proof. This technique could be extended to some range $p=2\pm c(n)$, but it does not seem to work for general $p>\frac{2n}{n+2}$. 
For this reason we state the boundary Calder\'on \& Zygmund estimates only for the case $p=2$. In this case \eqref{assump-b} simplifies to 
\begin{align}\label{assump-b-2}
\left\{
\begin{array}{c}
	|b(z,u,\xi)| + \big(1+|\xi|^2\big)^{\frac12}|\partial_\xi b(z,u,\xi)|
	\le
	L\,\big(1+|\xi|^2\big)\\[5pt]
	\langle\partial_\xi b(z,u,\xi)\xi_o,\xi_o\rangle
	\ge
	\nu\,|\xi_o|^2\\[5pt]
	|b(z,u,\xi) - b(z,u_o,\xi)|
	\le
	\theta\big(|u-u_o|\big)\big(1+|\xi|^2\big)
\end{array}
\right.
\end{align}
whenever $z\in\Omega_T$, $u,u_o\in\R^N$, $\xi,\xi_o\in\R^{Nn}$, for some parameters $0<\nu\le 1\le L$.
Then, we get the following result.

\begin{theorem}\label{thm:main-CZ2}
There exists $\epsilon=\epsilon(n,N,\nu,L)>0$ such that the following holds: whenever $q\in(2,2+\frac4n+\epsilon]$ and 
$$u\in C^0\big([0,T];L^{2} (\Omega,\R^N)\big)\cap L^2\big(0,T;W^{1,2}(\Omega,\R^N)\big)$$ 
is a weak solution to the parabolic Cauchy-Dirichlet problem \eqref{Cauchy-Dirichlet} with $p=2$, under the assumptions \eqref{def_asymp-CZ} -- \eqref{inhom} with \eqref{VMO-b} -- \eqref{assump-b-2} and 
\begin{equation*}
	b(z,u,\xi)\equiv b(z,\xi),
\end{equation*}
then we have
\begin{equation*}
	D u\in L^{q}\big(\Omega\times(\delta,T),\R^{Nn}\big)
	\quad\mbox{for any $\delta\in(0,T)$.}
\end{equation*}
Moreover, the quantitative estimate \eqref{CZ-est} holds with constants $c$ depending on $n,N,\nu,$ $L,q,\omega(\cdot),\V(\cdot),\partial\Omega$ and 
$R_o$ depending on 
$n, N,\nu,L,q, \|\omega\|_\infty, \V(\cdot),\partial\Omega$.
\end{theorem}

Without any regularity assumption on the solution $u$ it cannot be expected that the previous results hold with a vector field $b$ depending additionally on $u$, since then the mapping
$x\mapsto b(x,t,u(x,t),\xi)$ is only measurable. Therefore, in order to deal with the $u$-dependence we have to assume that the solution is continuous. Note that such an assumption can be verified for $p=2$ in the low-dimensional case $n=2$, see \cite{Campanato:1984, Boegelein:2012}.

\begin{theorem}\label{thm:main-CZ4}
There exists $\epsilon=\epsilon(n,N,\nu,L)>0$ such that the following holds: whenever $q\in(2,2+\frac4n+\epsilon]$ and 
$$u\in %C^0\big([0,T];L^{2} (\Omega,\R^N)\big)\cap 
L^2\big(0,T;W^{1,2}(\Omega,\R^N)\big) 
\cap C^{0}\big(\overline\Omega\times(0,T),\R^N\big)$$ 
is a weak solution to the parabolic Cauchy-Dirichlet problem \eqref{Cauchy-Dirichlet} with $p=2$, under the assumptions \eqref{def_asymp-CZ} -- \eqref{inhom} with \eqref{VMO-b} -- \eqref{assump-b-2},
then we have
\begin{equation*}
	D u\in L^{q}\big(\Omega\times(\delta,T),\R^{Nn}\big)
	\quad\mbox{for any $\delta\in(0,T)$.}
\end{equation*}
Moreover, the quantitative estimate \eqref{CZ-est}, respectively \eqref{CZ-est-int} holds with constants $c$ depending on $n,N,\nu,L,q,\omega(\cdot),\V(\cdot),\partial\Omega$ and 
$R_o$ on 
$n, N,\nu,L,q, \|\omega\|_\infty,\V(\cdot), \partial\Omega$.
\end{theorem}

Note that the preceding result can be used to obtain an improved global partial regularity result for parabolic systems in the case $n=2$, see \cite{Boegelein:2012}.

\subsection{The general setting}\label{sec:general-setting}
The preceding theorems are a consequence of the following abstract result together with suitable a priori estimates. Instead of assuming a certain regularity of the vector field $b$ we only assume that solutions to the associated parabolic system satisfy certain a priori estimates. This will be made precise in the following definition. For the notation we refer to Section \ref{sec:notation}.

\begin{definition}\label{def:chi-reg}
Let $\chi>p$.  A vector field $b\colon \Omega_T\times\R^N\times\R^{Nn}\to\R^N$satisfying \eqref{assump-b}
is called $\chi$-regular if there holds: 
Whenever $u_o\in\R^N$, $Q_{\rho,\lambda}(z_o)\subset \R^{n+1}$ with
$\rho\in(0,1]$, $\lambda \ge 1$ and $(x_o)_n\ge 0$ is a scaled parabolic cylinder, $\Phi\colon \overline{B_\rho(x_o)}\to \R^n$ a $C^1$-diffeomorphism with $\Phi(B_\rho^+(x_o))\subset\Omega$ and 
$$w\in C^0\big(\Lambda_{\rho,\lambda}(t_o);L^2(B_\rho^+(x_o),\R^N)\big) \cap L^p\big(\Lambda_{\rho,\lambda}(t_o);W^{1,p}(B_\rho^+(x_o),\R^N)\big)$$ 
is a weak solution to the frozen system
\begin{equation}\label{sys:chi-reg} 
	\left\{
	\begin{array}{cc}
	\partial_t w - \Div \big(b(\Phi(\cdot),t,u_o,Dw)\big)_{B_\rho^+(x_o)} = 0,
	&\mbox{in $Q_{\rho,\lambda}^+(z_o)$,}\\[7pt]
	w=0
	&
	\mbox{on $\Gamma_{\rho,\lambda}(z_o)$ if $\Gamma_{\rho,\lambda}(z_o)\not=\emptyset$,} 
	\end{array}
	\right.
\end{equation}
then for any $\kappa\ge 1$ the following improvement of integrability holds:
\begin{equation}\label{w-cond}
    \mint_{Q_{\rho,\lambda}(z_o)} |Dw|^p \dz
    \le
    \kappa\,\lambda^p
	\quad\Longrightarrow\quad
    \mint_{Q_{\rho/2,\lambda}(z_o)} |Dw|^\chi \dz
    \le
    H_b\,\lambda^\chi,
\end{equation}
for a constant $H_b(\kappa, \|\Phi\|_{C^1})$.
\end{definition}

Note that the $C^1$-diffeomorphism $\Phi$ in Definition \ref{def:chi-reg} is introduced in order to allow the flattening of the boundary procedure performed in Section \ref{sec:model-CZ} for cylinders intersecting the lateral boundary of $\Omega_T$. It can be omitted when only interior cylinders are considered. We could have given a definition for the interior and the boundary case separately. However, for our purposes it is more convenient to have a unified definition at hand, since we do not distinguish later on between boundary and interior cases.

The following abstract Calder\'on \& Zygmund theorem states that solutions to the original problem \eqref{Cauchy-Dirichlet} are almost as integrable as solutions to the more regular problem \eqref{sys:chi-reg}.

\begin{theorem}\label{thm:main-CZ}
Let $\chi>p$ and $q\in(p,\chi)$ and assume that
$$u\in C^0\big([0,T];L^{2} (\Omega,\R^N)\big)\cap L^p\big(0,T;W^{1,p}(\Omega,\R^N)\big)$$
is a weak solution to the partial Cauchy-Dirichlet problem \eqref{Cauchy-Dirichlet}, under the assumptions \eqref{def_asymp-CZ} -- \eqref{inhom} and \eqref{assump-b} -- \eqref{VMO-b-}.
Moreover, assume that the vector field $b$ is $\chi$-regular in the sense of Definition \ref{def:chi-reg} and if $b$ depends on $u$ we additionally assume that $u$ and $g$ are continuous in $\Omega_T$.
Then
\begin{equation*}
	D u\in L^{q}\big(\Omega\times(\delta,T),\R^{Nn}\big)
	\quad\mbox{for any $\delta\in(0,T)$.}
\end{equation*}
Moreover, the quantitative estimate \eqref{CZ-est}, respectively \eqref{CZ-est-int} hold with constants depending additionally on $H_b(\cdot)$.
\end{theorem}

Finally, we briefly comment on the strategy of proof.
The first ingredient is the intrinsic geometry of DiBenedetto \& Friedman mentioned already in the introduction. This technique was invented in \cite{DiBenedetto-Friedman:1984, DiBenedetto-Friedman:1985} to prove the $C^{1,\alpha}$-regularity of solutions to the parabolic $p$-Laplacian system. 
There, the brilliant idea of DiBenedetto \& Friedman was to introduce a system of parabolic cylinders different from the standard ones whose space-time scaling depends on the local behavior of the solution itself, and which, in a certain sense, rebalances the non-homogeneous scaling of the parabolic $p$-Laplacian system with respect to space and time. The strategy is to find so called \textit{intrinsic parabolic cylinders} of the form 
\begin{equation*}
	Q_{\rho,\lambda}(z_o)
	:=
	B_\rho(x_o)\times\big(t_o-\lambda^{2-p}\rho^2, t_o+\lambda^{2-p}\rho^2\big),
	\quad z_o=(x_o,t_o)
\end{equation*} 
in such a way that the scaling parameter $\lambda>0$ and the average of $|Du|^p$ over the cylinder are coupled by a condition of the type
\begin{equation*}
	\mint_{Q_{\rho,\lambda}(z_o)} |Du|^p \dz
	\approx
	\lambda^p.
\end{equation*} 
The delicate aspect within this coupling clearly relies in the fact that the value of the integral average must be comparable to the scaling factor $\lambda$ which itself is involved in the construction of its support. On such intrinsic cylinders the parabolic $p$-Laplacian system behaves in a certain sense like $\partial_t u = \lambda^{p-2}\Delta u$. Therefore, using cylinders of the type $Q_{\rho,\lambda}(z_o)$ allows to rebalance the occurring multiplicative factor $\lambda^{p-2}$ by rescaling $u$ in time with a factor $\lambda^{2-p}$.

The strategy now, is to consider the function $\tilde u:=u-g$ which has boundary values equal to zero and satisfies 
\begin{equation*}
	\partial_t \tilde u -
    \Div a(z,\tilde u+g,D\tilde u+Dg)
	=
	\Div(|F|^{p-2}F)+f-\partial_t g
	\quad\text{in } \Omega_T.
\end{equation*}
By a transformation argument we may assume that $\Omega_T$ is a half-cylinder.
Subsequently we cover some subset of $\Omega_T$ by intrinsic cylinders $Q_{\rho,\lambda}(z_o)$ satisfying
\begin{equation*}
	\mint_{Q_{\rho,\lambda}(z_o)\cap\Omega_T} |D\tilde u|^p \dz
	\approx
	\lambda^p.
\end{equation*} 
Such cylinders are constructed by an exit time argument.
On $Q_{\rho,\lambda}(z_o)\cap\Omega_T$ we compare $\tilde u$ to the solution $v$ of the homogeneous system
\begin{equation*}
	\partial_t v -
    \Div b(z,\tilde u+g,Dv)
	=
	0
	\quad\text{in } Q_{\rho,\lambda}(z_o)\cap\Omega_T
\end{equation*}
which has boundary values equal to $\tilde u$.
In a second step we compare $v$ to the solution $w$ of 
\begin{equation*}
	\partial_t w -
    \Div \mathcal B(t,Dw)
	=
	0
	\quad\text{in } Q_{\rho/2,\lambda}(z_o)\cap\Omega_T
\end{equation*}
which has boundary values equal to $v$. Thereby, $\mathcal B$ is defined as the average of $b$ with respect to the spatial direction and $u+g$ are replaced by their means. The advantage of the two step comparison technique is that we can use the higher integrability of $v$ to deal with the VMO-condition in the comparison estimate. 
From the a priori estimates for the function $w$ which is a solution to a more regular problem we infer estimates for super-level sets of $D\tilde u$. Since the maximal function is not compatible with the intrinsic geometry, we choose a maximal function free approach, which has its origin in \cite{Acerbi-Mingione:2007}. This finally yields the desired $L^q$-estimate.

\section{Preliminaries}

\subsection{Some notation}\label{sec:notation}
Throughout the paper we will generally write $x=(x_1,\dots,x_n)$ for a point in $\R^n$ and $z=(x,t)=(x_1,\dots,x_n,t)$ for a point in $\R^{n+1}$. By $B_\rho(x_o):=\{x\in\R^n:|x-x_o|<\rho\}$, respectively $B_\rho^+(x_o):=B_\rho(x_o)\cap \{x\in\R^n:x_n>0\}$ we denote the open ball, respectively upper part of the open ball in $\R^n$ with center $x_o\in\R^n$ and radius $\rho>0$. When considering $B_\rho^+(x_o)$ we do not necessarily assume $(x_o)_n=0$. Indeed, if $B_\rho(x_o)\subset\{x\in\R^n:x_n>0\}$ it can also happen that $B_\rho^+(x_o)\equiv B_\rho(x_o)$. Moreover, we write 
$$
	\Lambda_{\rho,\lambda}(t_o)
	:=
	\big(t_o-\lambda^{2-p}\rho^2, t_o+\lambda^{2-p}\rho^2\big)
$$
for the open interval around $t_o\in\R$ of length $2\lambda^{2-p}\rho^2$ with $\rho,\lambda>0$. As basic sets for our estimates we usually take cylinders. These are denoted by 
$$
	Q_{\rho,\lambda}(z_o)
	:=
	B_\rho(x_o)\times \Lambda_{\rho,\lambda}(t_o)
$$
and
$$
	Q_{\rho,\lambda}^+(z_o)
	:=
	B_\rho^+(x_o)\times \Lambda_{\rho,\lambda}(t_o),
$$
where $z_o=(x_o,t_o)\in\R^{n+1}$. As before, when considering $Q_{\rho,\lambda}^+(z_o)$ we do not necessarily assume $(x_o)_n=0$. 
For the hyperplane $x_n=0$ in $\R^{n+1}$ we write
$$
	\Gamma
	:=
	\big\{(x_1,\dots,x_{n-1},0,t)\in\R^{n+1}\big\}
$$
and
$$
	\Gamma_{\rho,\lambda}(z_o)
	:=
	Q_{\rho,\lambda}(z_o)\cap\Gamma
$$
for the flat part of the lateral boundary of $Q_{\rho,\lambda}^+(z_o)$. Note that it can happen that $\Gamma_{\rho,\lambda}(z_o)=\emptyset$. 
If $\lambda=1$ we use the shorter notations 
$$
	\Lambda_{\rho}(t_o)
	:=
	\Lambda_{\rho,1}(t_o),
	\quad
	Q_{\rho}(z_o)
	:=
	Q_{\rho,1}(z_o),
	\quad
	\Gamma_{\rho}(z_o)
	:=
	\Gamma_{\rho,1}(z_o)
$$
and if furthermore $z_o=0$ we write
$$
	B_{\rho}
	:=
	B_{\rho}(0),
	\quad
	\Lambda_{\rho}
	:=
	\Lambda_{\rho}(0),
	\quad
	Q_{\rho}
	:=
	Q_{\rho}(0),
	\quad
	\Gamma_{\rho}
	:=
	\Gamma_{\rho}(0).
$$

\subsection{Auxiliary tools}

The following lemma can be deduced from \cite[Lemma 2.1]{Acerbi-Fusco:1989}.
\begin{lemma}\label{lem:Fusco}
For every $\sigma\in(-1/2,0)$, $\mu\ge 0$ and $k\in\N$ we have
\begin{equation*}
    \int_0^1 \big(\mu^2 + |A + sB|^2\big)^\sigma\ds
    \le
    \frac{24}{2\sigma  +1}\,
    \big(\mu^2 + |A|^2 + |B|^2\big)^\sigma\,,
\end{equation*}
for any $A,B\in\R^{Nn}$, not both zero if $\mu=0$.
\hfill$\Box$
\end{lemma}

The next Lemma which is a consequence of \cite[Lemma 2.2]{Cupini-Fusco-Petti:1999} is a useful tool when working with $p$-growth problems; see also \cite[Lemma A.7]{Boegelein-habil}.
%****************************************************************
\begin{lemma}\label{lem:V}
Let $1<p<\infty$, $\mu\in[0,1]$ and $k\in\N$. Then there exists a constant $c\equiv c(k,p)$ such that

{\rm (i)}
for all $A,B\in\R^k$ there holds
\begin{equation*}
	\big(\mu^2 + |A|^2\big)^{\frac{p}{2}}
	\le
	c\,\big(\mu^2 + |B|^2\big)^{\frac{p}{2}} +
	c\,\big(\mu^2 + |A|^2 + |B|^2\big)^{\frac{p-2}{2}}|B-A|^2 ,
\end{equation*}

{\rm (ii)}
for all $A,B,C\in\R^k$ there holds
\begin{align*}
    \big(\mu^2 & + |A|^2 + |C|^2\big)^{\frac{p-2}{2}} |A-C|^2 \\
    &\le
    c\, \big(\mu^2 + |A|^2 + |B|^2\big)^{\frac{p-2}{2}} |A-B|^2 +
    c\, \big(\mu^2 + |B|^2 + |C|^2\big)^{\frac{p-2}{2}} |B-C|^2 .
\end{align*}
\end{lemma}

Next, we observe that the ellipticity assumption \eqref{assump-b}$_2$ on the vector field $b$ implies that $b$ is monotone, i.e.
\begin{align}\label{monotone}
	\big\langle b(z,u,\xi) - b(z,u,\xi_o), \xi-\xi_o\big\rangle
	%&=
%	\int_0^1
%	\partial_w a(\widetilde w + s(w - \widetilde w))
%	(w - \widetilde w)\cdot (w - \widetilde w) \ds \nonumber\\
%	&\ge
%	\nu
%	\int_0^1
%	\big(\mu^2+
%	|\widetilde w + s(w - \widetilde w)|^2\big)^{\frac{p-2}{2}}
%	|w - \widetilde w|^2 \,ds \nonumber\\
	&\ge
	\tfrac{\nu}{c(n,N,p)}\,
	\big(\mu^2 + |\xi|^2 + |\xi_o|^2\big)^{\frac{p-2}{2}}\,
	|\xi - \xi_o|^2 ,
\end{align}
holds for all $z\in \Omega_T$, $u\in\R^N$ and $\xi,\xi_o\in\R^{Nn}$. In particular, when $\xi_o=0$ we infer from \eqref{monotone} and Lemma \ref{lem:V} (i) (with $B=0$) that
\begin{equation}\label{monotone-2}
	\tfrac{\nu}{c(n,N,p)}\, |\xi|^p
	\le
	b(\xi)\cdot \xi + \mu^p.
\end{equation}

The following lemma is a parabolic analog of Sobolev's inequality. It is an immediate consequence of Gagliardo-Nirenberg's inequality, cf. \cite[Chapter I, Proposition 3.1]{DiBenedetto:book}

\begin{lemma}\label{lem:sob}
Let $u \in L^\infty(t_1,t_2;L^2(B))\cap L^p(t_1,t_2;W^{1,p}_0(B))$ with $p>1$, $t_1<t_2$ and 
$B:=B_\rho(x_o)$ be a ball in $\R^n$ with $0<\varrho\le1$. Then, we have
$u\in L^{\frac{p(n+2)}{n}}(B\times(t_1,t_2))$ and there exists a constant
$c=c(n,p)$ such that 
\begin{equation*}
	\int_{B\times(t_1,t_2)} |u|^\frac{p(n+2)}{n} \dz
	\le
	c\int_{B\times(t_1,t_2)} |Du|^p\dz\
	\bigg( \sup_{t\in(t_1,t_2)} \int_{B} |u(\cdot,t)|^2 \dx\bigg)
	^{\frac{p}{n}} .
\end{equation*}
\end{lemma}

Finally, we state an iteration lemma which is a standard tool in order to reabsorb certain terms from the right-hand side into the left, cf. \cite[Lemma 6.1]{Giusti:book}.
%*******************************************************************
%Iteration-Lemma
%*******************************************************************
\begin{lemma}\label{lem:it}
Let $\phi \colon [r,R]\to \R$ be a bounded non-negative function. Assume that
for $r\le s<t\le R$ there holds
\begin{equation*}
	\phi (t)
    \le
    \vartheta \phi (s)+\frac{A}{(t-s)^\alpha}+B,
\end{equation*}
where $A,B\ge 0$, $\alpha>0$ and $0\le\vartheta <1$. Then
\begin{equation*}
	\phi (r)
    \le
    c(\alpha ,\vartheta )\bigg[\frac{A}{(R-r)^\alpha}+B\bigg].
\end{equation*}
\end{lemma}

\subsection{Higher integrability}

We need the following Gehring's type higher integrability result for solutions to the more regular parabolic system associated to the vector field $b$. The interior result was first established by Kinnunen \& Lewis \cite{Kinnunen-Lewis:2000} and subsequently generalized in different directions, cf. \cite{Boegelein:2008, Boegelein-Parviainen:2010, Parviainen:2009}. Here, we need a version which covers both, the interior and the boundary situation.

\begin{lemma}\label{lem:hi}
There exists $\epsilon_o=\epsilon_o(n,N,p,\nu,L) >0$
such that the following holds: Whenever $Q_{\rho,\lambda}(z_o)\subset \R^{n+1}$ is a parabolic cylinder with $(x_o)_n\ge 0$, $\rho\in(0,1]$ and $\lambda\ge 1$ and 
$$
	v\in 
	C^0\big(\Lambda_{\rho,\lambda}(t_o);L^2(B_\rho^+(x_o),\R^N)\big)\cap 
	L^p\big(\Lambda_{\rho,\lambda}(t_o);W^{1,p}(B_\rho^+(x_o),\R^N)\big)
$$ 
is a weak solution
of the following partial Cauchy-Dirichlet problem
\begin{equation*}
\left\{
\begin{array}{cc}
	\partial_t v - \Div b(z,Dv) = 0
    &\mbox{in $Q_{\rho,\lambda}^+(z_o)$} \\[7pt]
    v=0
    &\mbox{on $\Gamma_{\rho,\lambda}(t_o)$ if\, $\Gamma_{\rho,\lambda}(t_o)\not=\emptyset$,}
\end{array}
\right.
\end{equation*}
where the vector field $b$ satisfies the following ellipticity and growth conditions
\begin{equation}\label{assump-b-hi}
	|b(z,\xi)|
    \le
    L(|\xi| + 1)^{p-1}, \quad
    \langle b(z,\xi),\xi\rangle
    \ge
    \nu |\xi|^p
    \qquad \forall\, z\in Q_{\rho,\lambda}^+(z_o), \xi\in\R^{Nn},
\end{equation}
then for any $\kappa\ge 1$ there exists a constant $c=c(n,N,p,\nu,L,\kappa)$ such that the bound
\begin{equation}\label{intrinsic-hi}
    \mint_{Q^+_{\rho,\lambda}(z_o)} |Dv|^{p} \dz
    \le
    \kappa\,\lambda^p ,
\end{equation}
implies
\begin{align*}
	\mint_{Q^+_{\rho/2,\lambda}(z_o)} |Dv|^{p(1+\epsilon_o)} \dz
	\le
	c\,\lambda^{p(1+\epsilon_o)} .
%\qquad\mbox{for any $\epsilon\in(0,\epsilon_o]$.}
\end{align*}
\end{lemma}

\begin{proof}
The qualitative higher integrability result, i.e. $Dv\in L^{p(1+\epsilon_o)}(Q^+_{\rho/2,\lambda}(z_o),\R^{Nn})$ directly follows from \cite{Kinnunen-Lewis:2000, Boegelein-Parviainen:2010}, but the quantitative estimates are stated there only on standard cylinders, i.e. for the case $\lambda=1$. 
The strategy to treat the general case is
to rescale to standard parabolic cylinders via a transformation in
time and then apply the known quantitative higher integrability estimates.
We define the rescaled solution
\begin{equation*}
    \tilde v(x,t)
    :=
    \lambda^{-1} \,
    v(x,\lambda^{2-p} t)
    \qquad\mbox{for $(x,t)\in Q_{\rho}^+(z_o)$}
\end{equation*}
and the rescaled coefficients
\begin{equation*}
    \tilde b(x,t,\xi)
    :=
    \lambda^{1-p}\, b(x,\lambda^{2-p} t,\lambda\xi)
    \qquad \mbox{for $(x,t)\in Q_{\rho}^+(z_o)$ and $\xi\in\R^{Nn}$.}
\end{equation*}
Note that $\tilde b$ still satisfies the growth and ellipticity assumptions \eqref{assump-b-hi} with the same ellipticity constant $\nu$ and upper bound $L$.
Then
\begin{equation*}
    \tilde v
    \in
    C^0\big(\Lambda_{\rho}(t_o);L^2(B_{\rho}^+(x_o),\R^N)\big)\cap
L^{p}\big(\Lambda_{\rho}(t_o);W^{1,p}(B_{\rho}^+(x_o),\R^N)\big)
\end{equation*}
is a weak solution to the parabolic system
\begin{equation*}
    \partial_t \tilde v -
    \Div \tilde b(\cdot, D\tilde v)
    =
    0
    \qquad\mbox{in $Q_{\rho}^+(z_o)$}
\end{equation*}
and moreover $\tilde v=0$ on $\Gamma_\rho(z_o)$ if $\Gamma_\rho(z_o)\not=\emptyset$. Therefore, \cite[Theorem 1]{Boegelein:2008} in the case $(x_o)_n\ge\rho$, respectively \cite[Theorem 2.2]{Boegelein-Parviainen:2010} in the case $0\le (x_o)_n<\rho$
yield the existence of an exponent $\epsilon_o=\epsilon_o(n,N,p,\nu,L) >0$ such that
\begin{align*}
	\mint_{Q_{\rho/2}^+(z_o)} |D\tilde v|^{p(1+\epsilon_o)} \dz
	&\le
	c \bigg(\mint_{Q_{\rho}^+(z_o)}
    |D\tilde v|^{p} \dz\bigg)^{1+\epsilon_o d} +
    c 
\end{align*}
holds for a constant $c=c(n,N,p,\nu,L)$.
Scaling back from $\tilde v$ to $v$ the preceding inequality turns into
\begin{align*}
	\mint_{Q_{\rho/2,\lambda}^+(z_o)} |Dv|^{p(1+\epsilon_o)} \dz
	&=
	\lambda^{p(1+\epsilon_o)}
    \mint_{Q_{\rho/2}^+(z_o)} |D\tilde v|^{p(1+\epsilon_o)} \dz \nn\\
	&\le
	c\, \lambda^{p(1+\epsilon_o)}\bigg[
    \bigg(\mint_{Q_{\rho}^+(z_o)}
    |D\tilde v|^{p} \dz\bigg)^{1+\epsilon_o d}  + 1\bigg] \nn\\
	&=
	c\, \lambda^{\epsilon_o p(1-d)}
    \bigg(\mint_{Q_{\rho,\lambda}^+(z_o)}
    |Dv|^{p} \dz\bigg)^{1+\epsilon_o d}  +
    c\,\lambda^{p(1+\epsilon_o)} \\
    &\le
    c\,\lambda^{p(1+\epsilon_o)} ,
\end{align*}
where in the last line we used hypothesis \eqref{intrinsic-hi}. This proves the assertion of the lemma.
\end{proof}

\subsection{A priori gradient estimates}\label{sec:apriori-CZ}
Here, we provide gradient estimates for solutions to the more regular parabolic systems. They will be used later to ensure that the vector field $b$ is $\chi$-regular in the sense of Definition \ref{def:chi-reg} with some $\chi>p$.
On a cylinder $Q_{\rho,\lambda}(z_o)\subset\R^{n+1}$ with $(x_o)_n\ge 0$, $\rho,\lambda\ge 1$ we consider Cauchy-Dirichlet problems of the type 
\begin{equation}\label{sys-apriori-1}
	\left\{
	\begin{array}{cc}
    \partial_t w - \Div b(t,Dw) = 0
    &\mbox{in $Q_{\rho,\lambda}^+(z_o)$} \\[5pt]
    w= 0
    &\!\! \mbox{on $\Gamma_{\rho,\lambda}(z_o)$ if $\Gamma_{\rho,\lambda}(z_o)\not=\emptyset$.}
    \end{array}
    \right.
\end{equation}
First, we consider a class of parabolic systems for which it is known that the spatial gradient of the solution is bounded. This result can be deduced from \cite{DiBenedetto:book}, Chapter VIII, Theorem 5.1 and Theorem 5.2' by a reflection argument. The argumentation is similar to 
\cite[Lemma 3.11]{Boegelein:2012-Lipschitz}.
Although the proof in \cite{DiBenedetto:book} is only for the $p$-Laplacian system, i.e. the case $c\equiv 1$, it can be adapted to our situation with minor changes.

\begin{lemma}\label{lem:sup-bound}
Let $Q_{\rho,\lambda}(z_o)\subset\R^{n+1}$ be a parabolic cylinder with $(x_o)_n\ge 0$, $\rho\in(0,1]$ and $\lambda\ge 1$. Further, suppose that 
$$w\in C^0\big(\Lambda_{\rho,\lambda}(t_o);L^2(B_\rho^+(x_o),\R^N)\big)\cap L^p\big(\Lambda_{\rho,\lambda}(t_o);W^{1,p}(B_\rho^+(x_o),\R^N)\big)$$ 
is a weak solution to \eqref{sys-apriori-1} in $Q_{\rho,\lambda}^+(z_o)$ with 
\begin{equation}\label{bb-1}
    b(t,\xi)
    \equiv
    c(t) |\xi|^{p-2}\xi
    \qquad\mbox{for $t\in \Lambda_{\rho,\lambda}(t_o)$ and $\xi\in\R^{Nn}$}
\end{equation}
and $c\colon \Lambda_{\rho,\lambda}(t_o)\to[\nu,L]$ measurable.
Then for any $c_\ast\ge 1$ there exists a constant $H(n,N,$ $p,\nu,L,c_\ast)\ge 1$ such that the bound
\begin{equation*}
    \mint_{Q_{\rho,\lambda}^+(z_o)} |Dw|^p \dz
    \le
    c_\ast\,\lambda^p
\end{equation*}
implies the following quantitative $L^\infty$-estimate
\begin{equation*}
    \sup_{Q_{\rho/2,\lambda}^+(z_o)} |Dw|
    \le
    H\,\lambda.
\end{equation*}
\end{lemma}

%\begin{proof}
%If $Q_{\rho/2,\lambda}(z_o)$ is an interior cylinder, i.e. $Q_{\rho/2,\lambda}(z_o)\subset Q_1^+$, the claim immediately follows 
%from \cite{DiBenedetto:book}, Theorem 5.1 in the case $p\ge 2$, respectively Theorem 5.2' in the case $p<2$
%applied to $w\circ \Psi_o^{-1}$.
%
%In the case that $Q_{\rho/2,\lambda}(z_o)$ is not an interior cylinder, i.e. $Q_{\rho/2,\lambda}(z_o)\setminus Q_1^+\not=\emptyset$ we use a reflection argument. 
%We define 
%\begin{equation*}
%    \tilde w(x,t)
%    =
%    \begin{cases}
%    w(x,t) & \mbox{if $x\in B^+$} \\[3pt]
%    -w(-x,t) & \mbox{if $x\in B\setminus B^+$}.
%    \end{cases}
%\end{equation*}
%Then, $\tilde w$ is a weak solution of 
%\begin{equation*}
%   	\partial_t \tilde w - \boldsymbol{c}\, \Div\big[\big(|D\tilde w\,\Psi_o|^{p-2}D\tilde w\,\Psi_o\big)\,\Psi_o^t\big]
%	=0
%	\qquad \mbox{in $Q_{\rho,\lambda}(z_o)$}.
%\end{equation*}
%Due to \eqref{w-cond-sup} and the fact that $z_o\in Q_1^+$ we have 
%\begin{equation*}
%    \mint_{Q_{\rho,\lambda}(z_o)} |D\tilde w|^p \dz
%    \le
%    2\mint_{Q_{\rho,\lambda}^+(z_o)} |Dw|^p \dz 
%    \le
%    2\kappa\,\lambda^p.
%\end{equation*}
%Then, we can proceed as in the interior case to conclude that there exists a constant $H(n,N,p,\nu,L,\kappa, \Psi_o)\ge 1$ such that 
%\begin{equation*}
%    \sup_{Q_{\rho/2,\lambda}^+(z_o)} |Dw|
%    \le
%    \sup_{Q_{\rho/2,\lambda}(z_o)} |D\tilde w|
%    \le
%    H\,\lambda.
%\end{equation*}
%This concludes the proof of the Theorem.
%\end{proof}

Next, we consider parabolic systems of the type
\eqref{sys-apriori-1},
where $b\colon\Lambda_{\rho,\lambda}(t_o)\times\R^{Nn}\to\R^{Nn}$ and $\partial_\xi b\colon\Lambda_{\rho,\lambda}(t_o)\times\R^{Nn}\to\R^{Nn}$ are Carath\'eodory maps satisfying 
\begin{equation}\label{assump-b-reg}
\left\{
\begin{array}{c}
	|b(t,\xi)| + \big(\mu^2+|\xi|^2\big)^{\frac12}|\partial_\xi b(t,\xi)|
	\le
	L\,(\mu^2+|\xi|^2)^{\frac{p-1}{2}}\\[5pt]
	\langle\partial_\xi b(t,\xi)\xi_o,\xi_o\rangle
	\ge
	\nu\,\big(\mu^2+|\xi|^2\big)^{\frac{p-2}{2}} |\xi_o|^2
\end{array}
\right.
\end{equation}
whenever $t\in \Lambda_{\rho,\lambda}(t_o)$ and $\xi,\xi_o\in\R^{Nn}$ for some constants $0<\nu\le L$ and $\mu\in[0,1]$.
Then, it is not anymore expected that the spatial gradient is bounded. Nevertheless, in the interior situation we have the following quantifiable higher integrability of the spatial gradient. In the case $p\ge 2$ it can be found in \cite[Lemma 5.9]{Duzaar-Mingione-Steffen:2011}, while in the case $\frac{2n}{n+2}<p<2$ it was proved in \cite[Theorem 4.3]{Scheven:2010} (here one can argue as in the proof of Lemma \ref{lem:hi} to extend the result from standard to general parabolic cylinders).

\begin{lemma}\label{lem:chi-bound-p}
There exists $\epsilon=\epsilon(n,N,p,\nu,L)>0$ such that the following holds:
Let $Q_{\rho,\lambda}(z_o)\subset \Omega_T$ be a parabolic cylinder with $z_o\in \Omega_T$, $\rho\in(0,1]$ and  $\lambda\ge 1$. Further, suppose that 
$$w\in C^0\big([0,T];L^2(\Omega,\R^N)\big)\cap L^p\big(0,T;W^{1,p}(\Omega,\R^N)\big)$$ 
is a weak solution to 
\begin{equation*}
    \partial_t w - b(t,Dw) = 0
    \qquad\mbox{in $\Omega_T$} 
\end{equation*}
under the assumptions \eqref{assump-b-reg} with $(0,T)$ instead of $\Lambda_1$.
Then, for any $c_\ast\ge 1$ there exists a constant $H(n,N,p,\nu,L,c_\ast)\ge 1$ such that the bound
\begin{equation*}
    \mint_{Q_{\rho,\lambda}(z_o)} |Dw|^p \dz
    \le
    c_\ast\lambda^p
\end{equation*}
implies 
\begin{equation*}
    \mint_{Q_{\rho/2,\lambda}(z_o)} |Dw|^\chi \dz
    \le
    H^\chi\lambda^\chi,
    \quad\mbox{where } \chi:=p+\tfrac4n + \epsilon.
\end{equation*}
\end{lemma}

In the case $p=2$ the quantifiable higher integrability of the spatial gradient can be extended up to the boundary. The proof relies on an up-to-the-boundary higher differentiability result from \cite{Boegelein-Duzaar-Mingione:2010-boundary-II}. It was proved by a delicate interplay between difference quotients in space and time.

\begin{lemma}\label{lem:chi-bound-2}
There exists $\epsilon=\epsilon(n,\nu,L)>0$ such that the following holds:
Let $Q_{\rho}(z_o)\subset\R^{n+1}$ be a parabolic cylinder with $(x_o)_n\ge 0$ and $\rho\in(0,1]$. Further, suppose that
$$w\in C^0\big(\Lambda_{\rho}(t_o);L^2(B_\rho^+(x_o),\R^N)\big)\cap L^2\big(\Lambda_{\rho}(t_o);W^{1,2}(B_\rho^+(x_o),\R^N)\big)$$ 
is a weak solution to \eqref{sys-apriori-1} in $Q_{\rho}^+(z_o)$ under the assumptions \eqref{assump-b-reg} with $p=2$.
Then, there holds
\begin{equation}\label{w-chi}
    \mint_{Q_{\rho/2}^+(z_o)} |Dw|^\chi \dz
    \le
    c\bigg(\mint_{Q_{\rho}^+(z_o)} (1+|Dw|)^2 \dz \bigg)^\frac\chi2
    \quad\mbox{where } \chi:=2+\tfrac4n + \epsilon
\end{equation}
and $c=c(n,N,\nu,L)$. 
\end{lemma}

\begin{proof}
In the case $(x_o)_n\ge \rho$ the result follows from \cite[Lemma 5.4 and 5.8]{Duzaar-Mingione-Steffen:2011} applied with $p=2$. We next consider the case where $z_o\in\Gamma$. From 
\cite[Lemma 4.11]{Boegelein-Duzaar-Mingione:2010-boundary-II} we infer that $Dw\in L^{2+\frac4n}(Q_{7\rho/8}^+(z_o),\R^{Nn})$ together with the estimate
\begin{equation}\label{w-2+}
    \mint_{Q_{7\rho/8}^+(z_o)} |Dw|^{2+\frac4n} \dz
    \le
    c(n,N,\nu,L)\bigg(\mint_{Q_{\rho}^+(z_o)} (1+|Dw|)^2 \dz \bigg)^{1+\frac2n}.
\end{equation}
Note that in the preceding estimate we slightly changed the radii, i.e. we replaced the radius $\rho/2$ from \cite[Lemma 4.11]{Boegelein-Duzaar-Mingione:2010-boundary-II} by $7\rho/8$ which is possible by a different choice of the involved cut-off functions. The same applies to the next two statements.
From \cite[Proposition 4.15]{Boegelein-Duzaar-Mingione:2010-boundary-II} we know that 
there exists $\sigma=\sigma(n,\nu,L)\in(2,2+\frac4n]$ such that $D^2w\in L^\sigma(Q_{3\rho/4}^+(z_o),\R^{Nn^2})$ and 
\begin{equation}\label{w-sigma}
    \int_{Q_{3\rho/4}^+(z_o)} |D^2w|^{\sigma} \dx
    \le
    c(n,N,\nu,L)\,\rho^{-\sigma}\int_{Q_{7\rho/8}^+(z_o)} (1+|Dw|)^\sigma \dz .
\end{equation}
Finally, from \cite[Proposition 4.3 and 4.10]{Boegelein-Duzaar-Mingione:2010-boundary-II} we have the following estimate
\begin{equation}\label{w-infty-2}
    \sup_{t\in\Lambda_{3\rho/4}(t_o)}\int_{B_{3\rho/4}^+(x_o)} |Dw(\cdot,t)|^{2} \dx
    \le
    c(n,N,\nu,L)\,\rho^{-2}\int_{Q_{\rho}^+(z_o)} (1+|Dw|)^2 \dz .
\end{equation}
Now, we choose a cut-off function $\eta\in C_0^\infty(B_{3\rho/4}(x_o),[0,1])$ with $\eta\equiv 1$ on $B_{\rho/2}(x_o)$ and $|D\eta|\le 8/\rho$. From Lemma \ref{lem:sob} we infer that
\begin{align*}
	& \mint_{Q_{3\rho/4}^+(z_o)} |\eta Dw|^{\frac{\sigma(n+2)}{n}} \dz \\
	&\phantom{mm}\le
	c(n,\sigma)\mint_{Q_{3\rho/4}^+(z_o)} |D(\eta Dw)|^{\sigma} \dz
	\bigg(\sup_{t\in\Lambda_{3\rho/4}(t_o)} \int_{B_{3\rho/4}^+(x_o)} |Dw(\cdot,t)|^2\dx \bigg)^{\frac{\sigma}{n}}.
\end{align*}
The first integral on the right-hand side is now further estimated with the help of \eqref{w-sigma}, H\"older's inequality and \eqref{w-2+} as follows:
\begin{align*}
	\mint_{Q_{3\rho/4}^+(z_o)} |D(\eta Dw)|^{\sigma} \dz
	&\le
	c \mint_{Q_{3\rho/4}^+(z_o)} |D^2w|^{\sigma} + \Big|\frac{Dw}{\rho}\Big|^\sigma\dz \\
	&\le
	c\,\rho^{-\sigma} \mint_{Q_{7\rho/8}^+(z_o)} (1+|Dw|)^{\sigma} \dz \\
	&\le
	c\,\rho^{-\sigma} \bigg(\mint_{Q_{7\rho/8}^+(z_o)} (1+|Dw|)^{2+\frac4n} \dz\bigg)^{\frac{n\sigma}{2(n+2)}} \\
	&\le
	c\,\rho^{-\sigma} \bigg(\mint_{Q_{\rho}^+(z_o)} (1+|Dw|)^{2} \dz\bigg)^{\frac{\sigma}{2}},
\end{align*}
where $c=c(n,N,\nu,L)$. Inserting this and \eqref{w-infty-2} above and taking into account that $\eta\equiv 1$ on $B_{\rho/2}(x_o)$ yields \eqref{w-chi} with the choice $\epsilon=(\sigma-2)(1+\frac2n)>0$. This proves the assertion in the case $z_o\in\Gamma$. 

Therefore it remains to consider the case where $0< (x_o)_n<\rho$.
Here, we cover $Q_{\rho/2}^+(z_o)$ by $M=M(n)$ cylinders $Q_{\rho/4}^+(z_i)$, $i=1,\dots,M$ such that $Q_{\rho/2}^+(z_i)\subset Q_{\rho}^+(z_o)$ and for any $i=1,\dots,M$ there holds either $(x_i)_n\ge \rho/4$ or $z_i\in\Gamma$. Then, by the preceding argumentation we can apply \eqref{w-chi} to $Q_{\rho/4}^+(z_i)$ for any $i=1,\dots,M$. Summing the resulting inequalities over $i=1,\dots,M$ yields
\begin{align*}
    \int_{Q_{\rho/2}^+(z_o)} |Dw|^\chi \dz
    &\le
    \sum_{i=1}^M \int_{Q_{\rho/4}^+(z_i)} |Dw|^\chi \dz \\
    &\le
    c\,\rho^{(n+2)(1-\frac{\chi}{2})}
    \sum_{i=1}^M \bigg(\int_{Q_{\rho/2}^+(z_i)} (1+|Dw|)^2 \dz \bigg)^\frac\chi2 \\
    &\le
    c\,\rho^{(n+2)(1-\frac{\chi}{2})}
    \bigg(\int_{Q_{\rho}^+(z_o)} (1+|Dw|)^2 \dz \bigg)^\frac\chi2.
\end{align*}
Taking mean values we obtain \eqref{w-chi} also in the case $0< (x_o)_n<\rho$ and this finishes the proof of the lemma. 
\end{proof}

\section{Proof of the Calder\'on \& Zygmund estimates}
Here, we first explain how the results from Sections \ref{sec:asym-plap} and \ref{sec:asym-gen} can be deduced from Theorem \ref{thm:main-CZ}. Therefore, we only have to ensure that the more regular vector field $b$ is $\chi$-regular in the sense of Definition \ref{def:chi-reg} with the appropriate value of $\chi$ and this is a consequence of the a priori estimates from Section \ref{sec:apriori-CZ}.
More precisely, in the setting of Theorem \ref{thm:main-CZ1} we can apply Lemma \ref{lem:sup-bound} to deduce that $b$ is $\infty$-regular -- which means that $b$ is $\chi$-regular for any $\chi<\infty$ -- 
while in the setting of 
Theorem \ref{thm:main-CZ3} we infer from Lemma \ref{lem:chi-bound-p} that $b$ is $p+\frac4n+\epsilon$-regular.
Finally, in the setting of Theorems \ref{thm:main-CZ2} and \ref{thm:main-CZ4} we apply Lemma \ref{lem:chi-bound-2} to find that $b$ is $2+\frac4n+\epsilon$-regular. 
This justifies the application of Theorem \ref{thm:main-CZ} with any $\chi<\infty$, respectively with $\chi=p+\frac4n+\epsilon$ or $\chi=2+\frac4n+\epsilon$ and therefore yields the results stated in Sections \ref{sec:asym-plap} and \ref{sec:asym-gen}.

The rest of the chapter is therefore devoted to the proof of the Calder\'on \& Zygmund theory in the most general setting stated in Theorem \ref{thm:main-CZ}.
The proof will be divided into several steps.

\subsection{Transformation to the model situation}\label{sec:model-CZ}

Since the asserted quantitative gradient bounds are of local nature we can locally transform the problem to a model situation on a cylinder intersected with the half-space and for boundary values $\equiv 0$ on the flat lateral boundary portion. The strategy is as follows.
We fix a cylinder $Q_R(z_o)\subset \R^n\times(0,T)$ such that $z_o=(x_o,t_o)\in\overline\Omega\times(0,T)$.
If $Q_{2R}(z_o)\subset \Omega_T$ we are in the interior situation and therefore can omit this step. Otherwise, if $Q_{2R}(z_o)\setminus \Omega_T\not=\emptyset$ we find $\tilde z=(\tilde x,\tilde t)\in \partial\Omega\times(0,T)\cap Q_{2R}(z_o)$. Without loss of generality we can assume that $\tilde z=0$
and that the inward pointing unit normal to $\partial\Omega$ in $\tilde x$ is $\nu_{\partial\Omega}(\tilde x) = e_n$. Note that $Q_{2R}(z_o)\subset Q_{4R}$.
Then, if $R>0$ is sufficiently small, we flatten the boundary $B_{2R}(x_o)\cap\partial\Omega$ by a
$C^{1}$-diffeomorphism  $\Phi$ as for instance constructed in \cite[Section 4.1]{Boegelein:2012-Lipschitz}, satisfying $\Phi(B_{2R}(x_o)\cap\partial\Omega)\subset B_{2R}(x_1) \cap \{x\in\R^n: x_n=0\}$ with $x_1:=\Phi(x_o)$ and 
\begin{equation*}
	B_{\rho/\sqrt{2}}^+(x_1)
	\subset
	\Phi(\Omega\cap B_\rho(x_o))
	\subset
	B_{\sqrt{2}\rho}^+(x_1)
	\qquad\mbox{for any $\rho\le \sqrt{8}R$}
\end{equation*}
and 
\begin{equation*}
	\Omega\cap B_{\rho/\sqrt{2}}(x_o)
	\subset
	\Phi^{-1}(B_\rho^+(x_1))
	\subset
	\Omega\cap B_{\sqrt{2}\rho}(x_o)
	\qquad\mbox{for any $\rho\le \sqrt{8}R$}
\end{equation*}
and $\det D\Phi=1=\det D\Phi^{-1}$. Next, we define the transformed maps
\begin{equation*}
	\hat g(y,t)
	:=
	g\big(\Phi^{-1}(y),t\big) ,\qquad 
	\hat F(y,t)
	:=
	F\big(\Phi^{-1}(y),t\big)
\end{equation*}
and
\begin{equation*}
	\hat f(y,t)
	:=
	f\big(\Phi^{-1}(y),t\big) + \partial_t\hat g(y,t)
\end{equation*}
and
$$
    v(y,t)
    :=
    u\big(\Phi^{-1}(y),t\big) - \hat g(y,t)
$$
for $(y,t)\in Q_{2R}^+(z_1)$, where $z_1:= (x_1,t_o)\equiv (\Phi(x_o),t_o)$.
Then, it can be shown that $v$ is a weak solution to the following Cauchy-Dirichlet problem
\begin{align*}
	\left\{
	\begin{array}{cc}
	\partial_t v -
    \Div \big[\hat a\big(z,v+\hat g,(Dv+D\hat g)\Psi\big)\Psi^t\big]
	=
	\Div\big[|\hat F|^{p-2}\hat F\,\Psi^t\big] + \hat f
	&
	\text{in } Q_{2R}^+(z_1)\\[7pt]
	v=0
	&
	\text{on } \Gamma_{2R}(z_1)
	\end{array}
	\right.
\end{align*}
where 
\begin{equation*}
    \Psi(y)
    :=
	D\Phi\big(\Phi^{-1}(y)\big)
\end{equation*}
and the vector-field $\hat a$ is defined by
\begin{equation*}
	\hat a(y,t,u,\xi)
	\equiv
	a\big(\Phi^{-1}(y),t, u,\xi\big).
\end{equation*}
The new vector field $\hat a$ is asymptotically related to 
\begin{equation*}
	\hat b(y,t,u,\xi)
	\equiv
	b\big(\Phi^{-1}(y),t, u,\xi\big)
\end{equation*}
in the sense of \eqref{def_asymp-CZ}, respectively Remark \ref{rem:asymp}, where $\Omega_T$ has to be replaced by $Q_{2R}^+(z_1)$.
From the assumptions \eqref{assump-b} -- \eqref{VMO-b-} on the vector-field $b\colon\Omega_T\times\R^N\times\R^{Nn}\to\R^N$, we infer that $\hat b\colon Q_{2R}^+(z_1)\times\R^N\times\R^{Nn}\to\R^N$
fulfills similar hypotheses with $Q_{2R}^+(z_1)$ instead of $\Omega_T$ after changing the appearing structure constants suitably.
More precisely, the new growth constant $\hat L$ then is of the form
$L\cdot c(\partial\Omega)$, while the new ellipticity constant
$\hat\nu$ is of the form $\nu/ c(\partial\Omega)$, where
the constant $c(\partial\Omega)$ is strictly larger then $0$.
Now, it is easy to verify that $Dv\in L^q(Q_{R}^+(z_1),\R^{Nn})$ if and only if
$Du\in L^q(\Phi^{-1}(Q_{R}^+(z_1)),\R^{Nn})$. Moreover, the quantitative estimate for $Du$ on $Q_R(z_o)\cap\Omega_T$ can be deduced from the one for $Dv$ together with a covering argument.
Therefore, it suffices to prove Theorem \ref{thm:main-CZ} in the model situation 
\begin{equation}\label{system-lat-CZ}
	\left\{
	\begin{array}{cl}
	\!\!\partial_t u -
    \Div \big[a\big(z,u{+}g,(Du{+}Dg)\Psi\big)\Psi^t\big]
	=
	\Div\big[|F|^{p-2}F \Psi^t\big] {+} f
	&
	\mbox{in $Q_{2R}^+(z_o)$} \\[7pt]
	u=0
	&
	\mbox{on $\Gamma_{2R}(z_o)$}
	\end{array}
	\right.
\end{equation}
for some cylinder $Q_{2R}(z_o)\subset\R^{n+1}$ with $(x_o)_n\ge 0$ and $0\in Q_{2R}^+(z_o)$. In the case $\Gamma_{2R}(z_o)\equiv Q_{2R}(z_o)\cap \Gamma=\emptyset$ the boundary condition \eqref{system-lat-CZ}$_2$ can obviously be omitted. The involved functions are assumed to satisfy 
\begin{equation}\label{Psi-CZ}
	\Psi\in C^{0}\big(B_{2R}^+(x_o)\cup D_{2R}(x_o),\R^n\big), 
	\ \Psi^{-1}\in L^\infty(B_{2R}^+(x_o),\R^n),\ \Psi(0)=\mathbb I_{n\times n}
\end{equation}
and 
\begin{equation}\label{inhom-lat}
	F,Dg\in L^q\big(Q_{2R}^+(z_o),\R^{Nn}\big),\quad
	f\in L^{\frac{q(n+2)}{p(n+2)-n}}\big(Q_{2R}^+(z_o),\R^N\big).
\end{equation}
Then, Theorem \ref{thm:main-CZ} is equivalent to the following

\begin{proposition}\label{prop:main-CZ}
There exists $R_o>0$ such that the following holds: Let $\chi>p$, $q\in(p,\chi)$ and $R\in(0,R_o]$ and assume that
$$
	u\in 
	C^0\big(\Lambda_{2R}(t_o);L^{2} (B_{2R}(x_o)^+,\R^N)\big)\cap 
	L^p\big(\Lambda_{2R}(t_o);W^{1,p}(B_{2R}(x_o)^+,\R^N)\big)
$$
is a weak solution to the partial Cauchy-Dirichlet problem \eqref{system-lat-CZ}, under the assumptions \eqref{def_asymp-CZ}, \eqref{assump-b} -- \eqref{VMO-b-}, \eqref{Psi-CZ} and \eqref{inhom-lat}.
Moreover, assume that the vector field $b$ is $\chi$-regular in the sense of Definition \ref{def:chi-reg} and if $b$ depends on $u$ we additionally assume that $u$ and $g$ are continuous in $Q_{2R}^+(z_o)\cup \Gamma_{2R}(z_o)$.
Then
\begin{equation*}
	D u\in L^{q}\big(Q_{R}^+(z_o),\R^{Nn}\big)
\end{equation*}
and moreover there holds
\begin{align*}
	\mint_{Q_{R}^+(z_o)} |D u|^{q} & \dz 
	\le
	c \Bigg[\bigg(\mint_{Q_{2R}^+(z_o)}|D u|^p\dz \bigg)^{\frac{q}{p}} +
	\mint_{Q_{2R}^+(z_o)}
	(|Dg| + |F|)^q \dz \nn\\
	&\phantom{mmmm}+
	R^{-(n+2)} \bigg(\int_{Q_{2R}^+(z_o)}
	|f|^{\frac{q(n+2)}{p(n+2)-n}} \dz 
	\bigg)^{1+\frac{p}{p(n+1)-n}} +1
	\Bigg]^{d_{CZ}},
\end{align*}
where $d_{CZ}$ is defined in \eqref{def-d}. Note that the constant $c$ depends on $n,N,p,\nu,L,\chi,q,\psi,$ $\omega(\cdot), H_b$ and $R_o$ depends on 
$n,N,p,\nu,L,\chi,q,\psi,$ $ \|\omega\|_\infty, H_b, \V(\cdot)$ and the modulus of continuity of $\Psi$ and if the vector field $b$ depends on $u$ also on $\theta(\cdot)$ and the moduli of continuity of $u$ and $g$.
Here, we have abbreviated
\begin{equation}\label{def-psi-CZ}
    \psi
    :=
    \|\Psi\|_{C^0(B_{2R}^+(x_o))} + \|\Psi^{-1}\|_{L^\infty(B_{2R}^+(x_o))}     
    \ge 1.
\end{equation}
\end{proposition}

%\begin{remark}\label{rem:main-CZ}\upshape
%In the case that $Q_{2R}(z_o)$ is an interior cylinder, we can omit the boundary transformation and directly consider problem \eqref{Cauchy-Dirichlet}. Then, the quantitative $L^q$-estimate is independent of the boundary data. More precisely, for a cylinder $Q_{2R}(z_o)\subset\Omega_T$ with $R\in(0,R_o]$ there holds
%\begin{align*}
%	\mint_{Q_{R}(z_o)} |D u|^{q}\dz
%	&\le
%	c\, \Bigg[\bigg(\mint_{Q_{2R}(z_o)}|D u|^p\dz \bigg)^{\frac{q}{p}} +
%	\mint_{Q_{2R}(z_o)} |F|^q \dz + 1\\
%	&\phantom{mmmm}+
%	R^{-(n+2)}
%	\bigg(\int_{Q_{2R}(z_o)} |f|^{\frac{q(n+2)}{p(n+2)-n}} \dz 
%	\bigg)^{1+\frac{p}{np+p-n}}
%	\Bigg]^{1+\frac{(d-1)(q-p)}{q}},
%\end{align*}
%where $c$ and $R_o$ are now independent of $\psi$.
%\end{remark}

With this respect, the rest of the chapter is devoted to the proof of Proposition \ref{prop:main-CZ}. 
The weak form of the partial Cauchy-Dirichlet problem \eqref{system-lat-CZ} reads as follows:
\begin{align}\label{system-lat-weak}
    \int_{Q_{2R}^+(z_o)} u\cdot\varphi_t  & -
    \big\langle a\big(z,u{+}g,(Du{+}Dg)\Psi\big),D\varphi\, \Psi\big\rangle \dz \nn\\
    &=
    \int_{Q_{2R}^+(z_o)} \langle |F|^{p-2}F,D\varphi\, \Psi\rangle - f\cdot \varphi \dz
    \quad\forall\, \varphi\in C_0^\infty(Q_{2R}^+(z_o),\R^N).
\end{align}

\subsection{Exit cylinders and covering}
We fix a cylinder $Q_{2R}(z_o)\subset\R^{n+1}$ with $(x_o)_n\ge 0$, $0\in Q_{2R}^+(z_o)$ and $R\in(0,R_o]$, where $R_o>0$ will be specified in the course of the proof, and suppose that the assumptions of Proposition \ref{prop:main-CZ} are in force.
With $d$ denoting the exponent from (\ref{def-d}) and $M\ge 1$ to be chosen later, we define $\lambda_o\ge 1$ by
\begin{equation}\label{def-lambda0}
    	\lambda_o
    	:=
    	\bigg[\mint_{Q_{2R}^+(z_o)}\big(|D u|^p + (M \G)^p + (M|f|)^{p_\#'} \boldsymbol{F} \big) \dz 
	\bigg]^{\frac{d}{p}} 
	\ge 1,
\end{equation}
where we have used the abbreviations
\begin{equation}\label{def-G}
    	\G(z)
    	:=
    	|Dg(z)| + |F(z)| + 1,
    	\quad
    	\boldsymbol{F}
	:=
    	\bigg(\int_{Q_{2R}^+(z_o)} (M|f|)^{p_\#'} \dz \bigg)^{\frac{p}{p(n+1)-n}}
\end{equation}
and 
\begin{equation*}
	p_\#:=\frac{p(n+2)}{n}
	\quad\mbox{and}\quad
	p'_\# 
	:= 
	\frac{p_\#}{p_\#-1}
	\equiv
	\frac{p(n+2)}{p(n+2)-n} \,.
\end{equation*}	
For fixed $R\le r_1<r_2\le 2R$ we consider the concentric parabolic
cylinders
\begin{equation*}
    Q_R^+(z_o)\subset Q_{r_1}^+(z_o)\subset Q_{r_2}^+(z_o)\subset Q_{2R}^+(z_o) .
\end{equation*}	
Then, for $\lambda>0$ parabolic cylinders of the type
$Q_{s,\lambda}^{+}(\mathfrak z_o) \equiv B_s^+(\mathfrak x_o)\times (\mathfrak t_o-\lambda^{2-p}s^2,\mathfrak t_o+\lambda^{2-p}s^2)$
with $\mathfrak z_o=(\mathfrak x_o,\mathfrak t_o)\in Q_{r_1}^+(z_o)$ and
$0<s\le \min\{\lambda^{\frac{p-2}{2}},1\} (r_2-r_1)/2$
are contained in $Q^+_{r_2}(z_o)$.
With the help of an exit time argument we now construct suitable intrinsic cylinders.
For $\lambda\ge 1$ we consider
the level set
\begin{equation}\label{def-levelset}
	E(\lambda,r_1)
	:=
	\big\{z\in Q_{r_1}^+(z_o):
    \mbox{$z$ is a Lebesgue point of $|Du|$ and $|Du(z)|>\lambda$}\big\}.
\end{equation}
By Lebesgue's differentiation theorem we know for any $\mathfrak z_o\in E(\lambda,r_1)$ that
\begin{align}\label{stop-lower}
	\lim_{s\downarrow 0} \bigg[
	\mint_{Q_{s,\lambda}^+(\mathfrak z_o)} \big(|D u|^p + (M \G)^p + (M|f|)^{p_\#'} \boldsymbol{F}\big) \dz 	
	\bigg] 
	\ge 
	|Du(\mathfrak z_o)|^p
	>
	\lambda^p.
\end{align}
In the following we consider radii $s$ such that
\begin{equation}\label{radius}
	\min\big\{\lambda^{\frac{p-2}{2}},1\big\}\,\frac{r_2-r_1}{2^{6}}
    \le
    s
    \le
    \min\big\{\lambda^{\frac{p-2}{2}},1\big\}\,\frac{r_2-r_1}{2}
\end{equation}
and values of $\lambda$ satisfying
\begin{equation}\label{def-B}
	\lambda
	>
	B\,\lambda_o,
	\qquad
	\text{where }\qquad
	B
    :=
    \Big(\frac{2^7 R}{r_2-r_1}\Big)^{\frac{d(n+2)}{p}} .
\end{equation}
By the definition of $\lambda_o$ from \eqref{def-lambda0} we then have
\begin{align*}
	\mint_{Q^+_{s,\lambda}(\mathfrak z_o)} &
	\big(|D u|^p + (M \G)^p + (M|f|)^{p_\#'} \boldsymbol{F}\big) \dz \\
	&\le
	\frac{|Q_{2R}^+(z_o)|}{|Q^+_{s,\lambda}(\mathfrak z_o)|}
	\mint_{Q_{2R}^+(z_o)} \big(|D u|^p + M^p\G^p + (M|f|)^{p_\#'} \boldsymbol{F}\big) \dz \\
	&\le
	\Big(\frac{2R}{s}\Big)^{n+2}
    \lambda^{p-2}\,(\lambda_o)^\frac{p}{d}.
\end{align*}
In the case $p\ge 2$ we have $\frac{p}{d} = 2$ and
$\min\{\lambda^{\frac{p-2}{2}},1\} = 1$. Using also the definition of $B$ from \eqref{def-B} we therefore find 
\begin{align*}
	\mint_{Q^+_{s,\lambda}(\mathfrak z_o)}
	\big(|D u|^p + (M \G)^p + (M|f|)^{p_\#'} \boldsymbol{F}\big) \dz 
	&\le
	\Big(\frac{2^7R}{r_2-r_1}\Big)^{n+2}
    \lambda^{p-2}\,\lambda_o^2 \nn\\
	&<
	\Big(\frac{2^7R}{r_2-r_1}\Big)^{n+2}
    B^{-2}\,\lambda^{p}
    =
    \lambda^p.
\end{align*}
On the other hand, in the case $p< 2$ we have
$\frac{p}{d} = [p(n+2)-2n]/2$ and
$\min\{\lambda^{\frac{p-2}{2}},1\} = \lambda^{\frac{p-2}{2}}$ which leads us to
\begin{align*}
	\mint_{Q^+_{s,\lambda}(\mathfrak z_o)}
	\big(|D u|^p + (M \G)^p + (M & |f|)^{p_\#'} \boldsymbol{F}\big) \dz 	
	\le
	\Big(\frac{2^7R\,\lambda^{\frac{2-p}{2}} }{r_2-r_1} \Big)^{n+2}
    \lambda^{p-2}\,(\lambda_o)^\frac{p}{d} \nn\\
	&<
	\Big(\frac{2^7R}{r_2-r_1} \Big)^{n+2}
    B^{-\frac{p}{d}}
    \lambda^{\frac{(2-p)(n+2)}{2}}
    \lambda^{p-2}\,\lambda^\frac{p}{d}
    =
    \lambda^p.
\end{align*}
Hence, in any case we have shown for $\mathfrak z_o\in Q_{r_1}^+(z_o)$ and $s$ and $\lambda$ chosen according to \eqref{radius} and \eqref{def-B} that
\begin{align}\label{stop-upper}
	\mint_{Q^+_{s,\lambda}(\mathfrak z_o)}
	\big(|D u|^p + (M \G)^p + (M|f|)^{p_\#'} \boldsymbol{F}\big) \dz 
	<
	\lambda^p.
\end{align}
From the preceding reasoning we conclude that \eqref{stop-lower} yields a radius for which the considered integral takes a value larger than $\lambda^p$, while  \eqref{stop-upper} states that the integral is smaller than $\lambda^p$ for any radius satisfying \eqref{radius}. Therefore, the continuity of the integral yields the existence of a maximal radius $\rho_{\mathfrak z_o}$ in between, i.e.
\begin{equation*}
    0 < \rho_{\mathfrak z_o}
    <
    \min\big\{\lambda^{\frac{p-2}{2}},1\big\}\,\frac{r_2-r_1}{2^{6}}
\end{equation*}
such that
\begin{align}\label{stop-radius}
	\mint_{Q^+_{\rho_{\mathfrak z_o},\lambda}(\mathfrak z_o)}
	\big(|D u|^p + M^p \G^p + (M|f|)^{p_\#'} \boldsymbol{F}\big) \dz 
	=
	\lambda^p
\end{align}
holds while
\begin{align}\label{stop-radius<}
	\mint_{Q^+_{s,\lambda}(\mathfrak z_o)}
	\big(|D u|^p + M^p \G^p + (M|f|)^{p_\#'} \boldsymbol{F}\big) \dz 
	<
	\lambda^p 
\end{align}
for any $s\in (\rho_{\mathfrak z_o},\,\min\{\lambda^{\frac{p-2}{2}},1\}\,\tfrac{r_2-r_1}{2}]$.
With this choice of $\rho_{\mathfrak z_o}$ we define concentric parabolic cylinders centered at $\mathfrak z_o$ as follows:
\begin{equation*}
	2^jQ^{+}_{\mathfrak z_o}
	:=
    	Q^+_{2^{j}\rho_{\mathfrak z_o},\lambda}(\mathfrak z_o)
	\subset
	Q^+_{r_2}(z_o)
	\qquad\text{for } j\in\{0,\dots,5\}.
\end{equation*}
At this point we remark that there are two possible cases included. Either $2^jQ_{\mathfrak z_o}$ intersects the hyperplane $\Gamma$, or $2^jQ_{\mathfrak z_o}$ is an interior cylinder, that is $2^jQ^+_{\mathfrak z_o}=2^jQ_{\mathfrak z_o}$.
From \eqref{stop-radius} and \eqref{stop-radius<} we conclude that
\begin{equation}\label{stop-Qj}
	\frac{\lambda^p}{2^{j(n+2)}}
	\le
	\mint_{2^jQ_{\mathfrak z_o}^+}
	\big(|D u|^p + M^p \G^p + (M|f|)^{p_\#'} \boldsymbol{F}\big) \dz 
	\le
	\lambda^p
	\quad\text{for } j\in\{0,\dots,5\}.
\end{equation}
Indeed, in the case $j=0$ \eqref{stop-Qj} is equivalent to (\ref{stop-radius}). In the case $j\in\{1,\dots,5\}$ the estimate from above in \eqref{stop-Qj} is exactly the one
from (\ref{stop-radius<}), while the estimate from below is a consequence of (\ref{stop-radius}), i.e.
\begin{align*}
	\mint_{2^jQ_{\mathfrak z_o}^+} &
	\big(|D u|^p + M^p \G^p + (M|f|)^{p_\#'} \boldsymbol{F}\big) \dz \\
	&\ge
	\frac{|Q_{\mathfrak z_o}^+|}{|2^jQ_{\mathfrak z_o}^+|}
	\mint_{Q_{\mathfrak z_o}^+}
	\big(|D u|^p + M^p \G^p + (M|f|)^{p_\#'} \boldsymbol{F}\big) \dz 
	=
	\frac{\lambda^p}{2^{j(n+2)}}.
\end{align*}
Thus we have shown that for every $\mathfrak z_o\in E(\lambda,r_1)$ there exists $0<\rho_{\mathfrak z_o}<\min\{\lambda^{\frac{p-2}{2}},1\}\frac{r_2-r_1}{2^{6}}$
such that on the parabolic cylinders $2^jQ_{\mathfrak z_o}^+$, $j\in\{0,\dots,5\}$ the estimate
(\ref{stop-Qj}) holds. Therefore, we are in position to apply Vitali's covering theorem. In this way we find an at most countable
family of disjoint parabolic cylinders $Q_i^+:=Q^+_{\rho_{\mathfrak z_i}, \lambda}(\mathfrak z_i)\subset Q_{r_2}^+(z_o)$
with center $\mathfrak z_i=(\mathfrak x_i,\mathfrak t_i)\in E(\lambda,r_1)$ such that the cylinders $2^3Q_{i}^+:=Q^+_{2^3\rho_{\mathfrak z_i}, \lambda}(\mathfrak z_i)\subset Q_{r_2}^+(z_o)$ cover
$E(\lambda ,r_1)$, i.e
\begin{equation}\label{levelset}
	E(\lambda,r_1)
	\subset
	\bigcup_{i=1}^\infty 2^3Q_i^+,
\end{equation}
and such that (\ref{stop-Qj}) is satisfied  with $\mathfrak z_i$ instead of $\mathfrak z_o$. For later use we introduce the following notations:
\begin{equation*}
    2^j Q_{i}^+
    :=
    2^j B_i^+\times 2^j\Lambda_i
    :=
    Q^+_{2^{j}\rho_{\mathfrak z_i},\lambda}(\mathfrak z_i)
    \subset
    Q_{r_2}^+(z_o)
\end{equation*}
with
\begin{equation*}
    2^j B_i^+
    :=
    B^+_{2^{j}\rho_{\mathfrak z_i}}(\mathfrak x_i)
    \quad\mbox{and}\quad
    2^j\Lambda_i
    :=
    \big(\mathfrak t_i-\lambda^{2-p}(2^{j}\rho_{\mathfrak z_i})^2, \mathfrak t_i+\lambda^{2-p}(2^{j}\rho_{\mathfrak z_i})^2\big).
\end{equation*}

\subsection{Comparison estimates}\label{sec:comparison-CZ}
Here, we shall proceed with a two step comparison technique. First, we compare the original solution $u$ to the solution $v$ of the homogeneous system associated to the asymptotic vector field $b$. In a second step, we compare $v$ to the solution $w$ of a frozen coefficient problem. The advantage of this two step procedure is that we can use the Gehring's type higher integrability result from Lemma \ref{lem:hi} for the first comparison function $v$ which does not apply to the original solution $u$ itself.
On the cylinder $2^5Q_i^+$ centered at $\mathfrak z_i=(\mathfrak x_i,\mathfrak t_i)$ we denote by
\begin{equation*}
    v_i\in
    C^0\big(2^5\Lambda_i;L^2(2^5 B_i^+,\R^N)\big)
    \cap
    L^p\big(2^5\Lambda_i;W^{1,p}(2^5 B_i^+,\R^N)\big)
\end{equation*}
the unique weak solution of the following parabolic Cauchy-Dirichlet problem:
\begin{equation}\label{system-compare-1}
	\left\{
	\begin{array}{cl}
	\partial_t v_i -
    \Div \big[b(z,u+g,Dv_i\,\Psi)\, \Psi^t\big]
	=
	0
	\quad &
	\text{in } 2^5Q_i^+\\[7pt]
	v_i=u
	\quad &
	\text{on } \partial_p 2^5Q_i^+.
	\end{array}
	\right.
\end{equation}
In the following we shall derive suitable energy and comparison estimates for $v_i$. Thereby, we shall argue again somewhat formal concerning the use of the time derivatives $\partial_tu$ and $\partial_tv_i$. The arguments can be made rigorous by a standard mollifying procedure as for instance by the use of Steklov averages.
Subtracting the weak form \eqref{system-lat-weak} of \eqref{system-lat-CZ} from the one of \eqref{system-compare-1} we get
\begin{align*}
    & - \int_{2^5Q_i^+} 
    (v_i -u)\cdot\partial_t\varphi \dz +
    \int_{2^5Q_i^+} \big\langle b(z,u+g,Dv_i\,\Psi),
    D\varphi\, \Psi\big\rangle \dz \nn\\
    &\phantom{mmm} =
    \int_{2^5Q_i^+} \big\langle a\big(z,u+g,(Du+Dg)\Psi\big),D\varphi\,\Psi\big\rangle -
    \big\langle|F|^{p-2}F, D\varphi\,\Psi\big\rangle + f\cdot\varphi \dz
\end{align*}
whenever $\varphi\in C^\infty_0(2^5Q_i^+,\R^N)$.
For $\tau\in 2^5\Lambda_i$ and $\epsilon>0$ such that $\tau+\epsilon\in 2^5\Lambda_i$ we define
\begin{equation}\label{chi-epsilon}
	\chi_\epsilon(t)
	:=
	\begin{cases}
		1 & \mbox{on $(-\infty, \tau]$,} \\[2pt]
		1-\frac1\epsilon(t-\tau) & \mbox{on $(\tau,\tau+\epsilon)$,} \\[2pt]
		0 & \mbox{on $[\tau+\epsilon,\infty)$.}
	\end{cases}
\end{equation}
In the preceding identity we then use the testing-function $\varphi_\epsilon=(v_i-u)\chi_\epsilon$.
Then, for a.e. $\tau\in 2^5\Lambda_i$ we get for the first integral on the left-hand side
\begin{align}\label{time-part}
    -&\int_{2^5Q_i^+} (v_i-u)\cdot\partial_t\varphi_\epsilon \dz \\
    &=
    -\int_{2^5Q_i^+}  \partial_t\big(|v_i-u|^2 \chi_\epsilon\big) - \partial_t(v_i-u) \cdot (v_i-u)\chi_\epsilon \dz \nn\\
    &=
    \frac12 \int_{2^5Q_i^+} \partial_t|v_i-u|^2 \chi_\epsilon \dz 
    =
    - \frac12 \int_{2^5Q_i^+} |v_i-u|^2 \partial_t\chi_\epsilon \dz \nn\\
    &=
    \frac1{2\epsilon} \int_\tau^{\tau+\epsilon} \int_{2^5B_i^+} |v_i-u|^2 \dx\dt 
    \to
    \frac1{2} \int_{2^5B_i^+} |(v_i-u)(\cdot,\tau)|^2 \dx
    \quad \mbox{as $\epsilon\downarrow 0$.} \nn
\end{align}
Hence in the limit $\epsilon\downarrow 0$ we obtain
\begin{align}\label{test-diff}
    \frac1{2} \int_{2^5B_i^+} |(v_i-u)(\cdot,\tau)|^2 & \dx +
    \int_{Q_\tau}\big\langle b(z,u+g ,Dv_i\,\Psi),
    D(v_i-u) \Psi\big\rangle \dz \\
    & \le
    \int_{Q_\tau} \big\langle a\big(z,u+g,(Du+Dg)\Psi\big),D(v_i-u)\Psi\big\rangle \dz \nn\\
    &\phantom{\le\ } -
    \int_{Q_\tau} \big\langle|F|^{p-2}F, D(v_i-u)\Psi\big\rangle - f\cdot(v_i-u) \dz,
    \nn
\end{align}
where we have abbreviated 
$Q_\tau := 2^5Q_i^+ \cap \{(x,t) \in \R^{n+1}: t\le\tau\}$.
Inequality \eqref{test-diff} will be used in the following in two different directions. First, we shall derive an energy bound for $Dv_i$. Therefore, we observe that \eqref{def_asymp-CZ-} and \eqref{assump-b}$_1$ imply
\begin{equation*}
    |a(z,u,\xi)|
    \le
    \big(\|\omega\|_\infty + L\big)\,(1+|\xi|)^{p-1}
    \quad\forall\,z\in Q_{2R}^+(z_o),u\in\R^N,\xi\in\R^{Nn}.
\end{equation*}
Using also \eqref{monotone-2}, recalling the definition of the constant $\psi$ from \eqref{def-psi-CZ}  and using Young's inequality we deduce from \eqref{test-diff} that for $\vartheta>0$ there holds
\begin{align*}
    \frac1{2} \int_{2^5B_i^+} & |(v_i-u)(\cdot,\tau)|^2 \dx +
    \frac{\nu}{c(n,N,p)\psi^p} \int_{Q_\tau} |Dv_i|^p \dz \\
    &\le
    \frac1{2} \int_{2^5B_i^+} |(v_i-u)(\cdot,\tau)|^2 \dx +
    \int_{Q_\tau}\big\langle b(z,u+g,Dv_i\,\Psi),
    Dv_i\, \Psi\big\rangle + \mu^p \dz \nn\\
    & \le
    L\int_{Q_\tau} \big(\mu^2 + |Dv_i\,\Psi|^2\big)^{\frac{p-1}{2}} |Du\, \Psi| + \mu^p\dz \\
    &\phantom{\le\ }+
    \big(\|\omega\|_\infty + L\big)
    \int_{Q_\tau} \big(1 + |(Du+Dg)\Psi|\big)^{p-1} |(Dv_i-Du)\Psi| \dz \nn\\
    &\phantom{\le\ } +
    \int_{Q_\tau} |F|^{p-1} |D(v_i-u)\Psi| + |f| |v_i-u| \dz \\
    & \le
    \vartheta\int_{Q_\tau} |Dv_i|^p \dz +
    c\int_{Q_\tau} \big(|Du|^p + |Dg|^p + |F|^p + 1\big) \dz \\
    &\phantom{\le\ }+
    c\int_{Q_\tau} |f| |v_i-u| \dz,
\end{align*}
where $c=c(n,N,p,L,\psi,\|\omega\|_\infty,1/\vartheta)$.
It remains to estimate the last integral on the right-hand side. To this aim we use H\"older's inequality, the Sobolev-type inequality from Lemma \ref{lem:sob} and Young's inequality to infer that
\begin{align*}
    \int_{Q_\tau}\! |f| |v_i-u| \dz
    &\le
    \bigg(\int_{Q_\tau} |v_i-u|^{p_\#} \dz \bigg)^{\frac{1}{p_\#}}
    \bigg(\int_{Q_\tau} |f|^{p_\#'} \dz \bigg)^{\frac{1}{p_\#'}} \\
    &\le
    c\bigg(\int_{Q_\tau}\! |Dv_i-Du|^{p} \dz \bigg)^{\frac{1}{p_\#}}
    \bigg(\sup_{t\in 2^5\Lambda_i} \int_{2^5B_i^+}\! |(v_i-u)(\cdot,t)|^{2} \dx \bigg)^{\frac{1}{n+2}} \\
    &\phantom{\le\ } \cdot
    \bigg(\int_{Q_\tau} |f|^{p_\#'} \dz \bigg)^{\frac{1}{p_\#'}} \\
    &\le
    \vartheta \int_{Q_\tau} \big(|Dv_i|^p + |Du|^{p}\big) \dz +
    \vartheta \sup_{t\in 2^5\Lambda_i} \int_{2^5B_i^+} |(v_i-u)(\cdot,t)|^{2} \dx  \\
    &\phantom{\le\ } +
    c(n,N,p) \bigg(\int_{Q_\tau} |f|^{p_\#'} \dz \bigg)^{\frac{p(n+2)-n}{p(n+1)-n}} .
\end{align*}
For the last integral of the right-hand side we use the definition of $\boldsymbol{F}$ from \eqref{def-G} to estimate
\begin{align*}
    \bigg(\int_{Q_\tau} |f|^{p_\#'} \dz \bigg)^{\frac{p(n+2)-n}{p(n+1)-n}} 
    \le
    \boldsymbol{F} \int_{Q_\tau} |f|^{p_\#'} \dz.
\end{align*}
Inserting this above and choosing $\vartheta=\nu/(4c(n,N,p)\psi^p)\le \frac14$ we can absorb the integral involving $|Dv_i|^p$ from the right-hand side into the left. This leads us to 
\begin{align*}
    \frac1{2} & \int_{2^5B_i^+} |(v_i-u)(\cdot,\tau)|^2 \dx +
    \frac{\nu}{2 c(n,N,p)\psi^p} \int_{Q_\tau} |Dv_i|^p \dz \\
    & \le
    \frac14 \sup_{t\in 2^5\Lambda_i} \int_{2^5B_i^+} |(v_i-u)(\cdot,t)|^{2} \dx \\
    &\phantom{\le\ } +
    c\int_{2^5Q_i^+} \big(|Du|^p + |Dg|^p + |F|^p + |f|^{p_\#'} \boldsymbol{F}  + 1\big)\dz,
\end{align*}
where $c=c(n,N,p,\nu,L,\psi,\|\omega\|_\infty)$.
Since the preceding inequality holds for a.e. $\tau\in 2^5\Lambda_i$ we can use it in two directions. In the second term on the left-hand side we let $\tau\uparrow t_i+\lambda^{2-p}(2^5\rho_{z_i})^2$ while in the first one we take the supremum over $\tau\in 2^5\Lambda_i$. This allows us to absorb the first term of the right-hand side into the left. 
Summing up the two resulting inequalities we obtain 
\begin{align*}
    \sup_{t\in 2^5\Lambda_i} & \int_{2^5B_i^+} |(v_i-u)(\cdot,t)|^{2} \dx +
    \int_{2^5Q_i^+} |Dv_i|^p \dz \nn \\
    &\le
    c\int_{2^5Q_i^+} \big(|Du|^p + |Dg|^p + |F|^p + |f|^{p_\#'} \boldsymbol{F} + 1\big) \dz.
\end{align*}
Using the intrinsic coupling from \eqref{stop-Qj} 
we obtain the following {\it energy estimate for $Dv_i$}:
\begin{align}\label{energy-v}
    \sup_{t\in 2^5\Lambda_i} \int_{2^5B_i^+} |(v_i-u)(\cdot,t)|^{2} \dx +
    \int_{2^5Q_i^+} |Dv_i|^p \dz 
    \le
    c\,\lambda^p\,|Q_i^+|
\end{align}
with a constant $c=c(n,N,p,\nu,L,\psi,\|\omega\|_\infty)$.

Next, we shall prove a comparison estimate for $Dv_i-Du$. Here, we start from inequality \eqref{test-diff}, rewritten in the following form:
\begin{align}\label{compare-1}
    &\int_{2^5Q_i^+}\big\langle b(z,u+g,Dv_i\,\Psi) - b(z,u+g,Du\,\Psi),
    D(v_i-u)\, \Psi\big\rangle \dz \\
    &\ \le
    \int_{2^5Q_i^+}
    \big\langle b\big(z,u+g,(Du+Dg)\Psi\big) - b\big(z,u+g,Du\,\Psi\big),
    D(v_i-u) \Psi\big\rangle \dz \nn\\
    &\ \phantom{\le\ }+
    \int_{Q^+}
    \big\langle a\big(z,u{+}g,(Du{+}Dg)\Psi\big) - b\big(z,u{+}g,(Du{+}Dg)\Psi\big)
,D(v_i{-}u)\Psi\big\rangle \dz \nn\\
    &\ \phantom{\le\ } -
    \int_{2^5Q_i^+} \big\langle|F|^{p-2}F, D(v_i-u)\Psi\big\rangle \dz +
    \int_{2^5Q_i^+} f\cdot (v_i-u) \dz \nn\\
    &\ =:
    \mbox{I} + \mbox{II} + \mbox{III} + \mbox{IV}.
    \nn
\end{align}
Here, we have taken into account that the first term on the right-hand side of \eqref{test-diff} is non-negative.
For the integral on the left-hand side of \eqref{compare-1}
by \eqref{monotone} we have the following lower bound:
\begin{align}\label{compare-r2}
    \int_{2^5Q_i^+} &
    \big\langle b(z ,u+g,Dv_i\, \Psi) - b(z,u+g,Du\,\Psi),
    D(v_i-u)\Psi\big\rangle \dz \\
    &\ge
    \frac{\nu}{c(n,N,p)\psi^p}\int_{2^5Q_i^+}
    \big(\mu^2 + |Dv_i|^2 + |Du|^2\big)^{\frac{p-2}{2}} |Dv_i-Du|^2 \dz .
    \nn
\end{align}
Next, we in turn estimate the terms I -- IV on the right-hand side of \eqref{compare-1}.
From \eqref{assump-b}$_1$, Lemma \ref{lem:Fusco}, H\"older's inequality, the energy estimate \eqref{energy-v}, \eqref{stop-Qj} and the facts that $\mu\le 1$ and $\lambda\ge 1$ we get
\begin{align*}
    \mbox{I}
    &=
    \int_{2^5Q_i^+}\int_0^1 
    \big\langle D_\xi b\big(z,u+g,Du\,\Psi + sDg\,\Psi\big) Dg\,\Psi,
    D(v_i-u) \Psi\big\rangle \ds\dz
    \\
    &\le
    L\int_{2^5Q_i^+}\int_0^1 
    \big(\mu^2 + |Du\,\Psi + sDg\,\Psi|^2\big)^{\frac{p-2}{2}} |Dg\,\Psi| |D(v_i-u) \Psi| \ds\dz
    \\
    &\le
    c(p)\,L \int_{2^5Q_i^+}
    \big(\mu^2+|Du\,\Psi|^2 + |Dg\,\Psi|^2\big)^{\frac{p-2}{2}}
    |Dg\,\Psi| |D(v_i-u)\,\Psi| \dz\\
    &\le
    c(p)\,L\,\psi^{p} \int_{2^5Q_i^+}
    \big(\mu^2+|Du|^2 + |Dg|^2\big)^{\frac{p-2}{2}} |Dg| |Dv_i-Du| \dz \\
    &\le
    c(p)\,L\,\psi^{p} \int_{2^5Q_i^+}
    |Dg|^{p-1} |Dv_i-Du| + 
    \chi_{p>2}\big(\mu^2+|Du|^2\big)^{\frac{p-2}{2}} |Dg| |Dv_i-Du| \dz \\
    &\le
    c(p)\,L\,\psi^{p} \int_{2^5Q_i^+}
    |Dg|^{p-1} \big(|Dv_i|+|Du|\big) +
    \big(\mu^2 + |Dv_i|^2+|Du|^2\big)^{\frac{p-1}{2}} |Dg| \dz \\
    &\le
    c(p)\,L\,\psi^{p}
    \bigg(\int_{2^5Q_i^+} \big(|Dv_i|^p+|Du|^p\big) \dz \bigg)^{\frac1p}
    \bigg(\int_{2^5Q_i^+} |Dg|^{p} \dz \bigg)^{1-\frac1p} \\
    &\phantom{\le\ } +
    c(p)\,L\,\psi^{p}
    \bigg(\int_{2^5Q_i^+}
    \big(\mu^p + |Dv_i|^p+|Du|^p\big) \dz\bigg)^{1-\frac1p}
    \bigg(\int_{2^5Q_i^+} |Dg|^p \dz \bigg)^{\frac1p} \\
    &\le
    c(n,N,p,\nu,L,\psi,\|\omega\|_\infty)\,
    \bigg[\frac{1}{M} + \frac{1}{M^{p-1}}\bigg]\,\lambda^p\,|Q_i^+|.
\end{align*}
Next, we come to the estimate of the term II. Here, we let $\delta\in(0,1]$ to be chosen later. Due to hypothesis
\eqref{omega-CZ} there exists $K_\delta>0$ depending on $\delta$ and $\omega$ such that $\omega(s)\le \delta$ for any $s\ge K_\delta-1$. In turn this implies 
\begin{equation*}
    \omega(s)(1+s)^{p-1}
    \le
    \delta\,(1+s)^{p-1} +
    \|\omega\|_\infty K_\delta^{p-1}
    \qquad \forall\, s\ge 0.
\end{equation*}
Hence, using \eqref{def_asymp-CZ-}, the preceding inequality and Young's inequality we find
\begin{align*}
    \mbox{II}
    &\le
    \int_{2^5Q_i^+}
    \big|b\big(z,u+g,(Du+Dg)\Psi\big) - a\big(z,u+g,(Du+Dg)\Psi\big)\big|
    |D(v_i-u)\,\Psi| \dz \nn\\
    &\le
    \int_{2^5Q_i^+}
    \omega\big(|(Du+Dg)\Psi|\big) \big(1+|(Du+Dg)\Psi|\big)^{p-1}
    |D(v_i-u)\,\Psi| \dz \nn\\
    &\le
    \psi
    \int_{2^5Q_i^+} \Big[\delta\, \big(1+|(Du+Dg)\Psi|\big)^{p-1} +
    \|\omega\|_\infty K_\delta^{p-1}\Big]
    |Dv_i-Du| \dz \nn\\
    &\le
    c(p)\, \psi^{p} \, \delta
    \int_{2^5Q_i^+} \big(1+|Du|^p+|Dg|^{p} + |Dv_i|^p\big) \dz \\
    &\phantom{\le\ } +
    c(p)\,\psi\, \|\omega\|_\infty K_\delta^{p-1} |Q_i^+|^{1-\frac1p}
    \bigg(\int_{2^5Q_i^+} \big(|Dv_i|^p + |Du|^p\big) \dz\bigg)^{\frac1p} .
\end{align*}
Using the energy estimate \eqref{energy-v} and \eqref{stop-Qj} we get
\begin{align*}
    \mbox{II}
    \le
    c\,
    \Big[\delta + \Big(\frac{K_\delta}{\lambda}\Big)^{p-1}\Big]\, \lambda^p \, |Q_i^+|
    \le
    c(n,N,p,\nu,L,\psi,\|\omega\|_\infty)\,
    \delta \, \lambda^p \, |Q_i^+|,
\end{align*}
provided
\begin{equation}\label{lambda-asym}
    \lambda
    \ge
    \frac{K_\delta}{\delta^{\frac{1}{p-1}}}\,.
\end{equation}
For the estimate of the term III we use H\"older's inequality, the energy estimate \eqref{energy-v} and \eqref{stop-Qj}, yielding that
\begin{align*}
    \mbox{III}
    &\le
    \psi\int_{2^5Q_i^+} |F|^{p-1} |Dv_i-Du| \dz \\
    &\le
    c(p)\,\psi \bigg(\int_{2^5Q_i^+} |F|^p \dz \bigg)^{1-\frac1p}
    \bigg(\int_{2^5Q_i^+} \big(|Du|^p + |Dv_i|^p\big) \dz \bigg)^{\frac1p} \\
    &\le
    c(n,N,p,\nu,L,\psi,\|\omega\|_\infty)\,
    \frac{\lambda^p}{M^{p-1}} \, |Q_i^+|.
\end{align*}
Finally, the term IV is estimated with the help of H\"older's inequality as follows:
\begin{align*}
    \mbox{IV}
    &\le
    \int_{2^5Q_i^+} |f| |v_i-u| \dz \\
    &\le
    \bigg(\int_{2^5Q_i^+} |v_i-u|^{\frac{p(n+2)}{n}} \dz \bigg)^{\frac{n}{p(n+2)}}
    \bigg(\int_{2^5Q_i^+} |f|^{\frac{p(n+2)}{p(n+2)-n}} \dz \bigg)^{\frac{p(n+2)-n}{p(n+2)}} 
    =:
    \mbox{IV}_1\cdot \mbox{IV}_2
\end{align*}
with the obvious meaning of IV$_1$ and IV$_2$.
Applying the Sobolev-type inequality from Lemma \ref{lem:sob}, the energy estimate \eqref{energy-v} and \eqref{stop-Qj} we obtain
\begin{align*}
    \mbox{IV}_1
    &\le
    c(n,p)\bigg(\int_{2^5Q_i^+} |Dv_i-Du|^{p} \dz \bigg)^{\frac{n}{p(n+2)}}
    \bigg(\sup_{t\in 2^5\Lambda_i} \int_{2^5B_i^+} |(v_i-u)(\cdot,t)|^{2} \dx 
    \bigg)^{\frac{1}{n+2}} \\
    &\le
    c\,(\lambda^p |Q_i^+|)^{\frac{n+p}{p(n+2)}} ,
\end{align*}
where $c=c(n,N,p,\nu,L,\psi,\|\omega\|_\infty)$.
Recalling that $p_\#'=p(n+2)/(p(n+2)-n)$, $2^5Q_i^+\subset Q_{2R}^+(z_o)$ and the definition of $\boldsymbol{F}$ and finally using \eqref{stop-Qj} we find that
\begin{align*}
    \mbox{IV}_2
    &=
    \frac1M 
    \bigg(\int_{2^5Q_i^+} (M|f|)^{\frac{p(n+2)}{p(n+2)-n}} \dz \bigg)^{\frac{p(n+2)-n}{p(n+2)}}
    \\
    &=
    \frac1M 
    \bigg(\int_{2^5Q_i^+} (M|f|)^{p_\#'} \dz \bigg)^{\frac{p(n+2)-n}{p(n+2)}} \\
    &=
    \frac1M 
    \bigg(\int_{2^5Q_i^+} (M|f|)^{p_\#'} \dz \bigg)^{\frac{p}{p(n+2)}} 
    \bigg(\int_{2^5Q_i^+} (M|f|)^{p_\#'} \dz \bigg)^{\frac{p(n+1)-n}{p(n+2)}}\\
    &\le
    \frac1M \, 
    \boldsymbol{F}^{\frac{p(n+1)-n}{p(n+2)}} 
    \bigg(\int_{2^5Q_i^+} (M|f|)^{p_\#'} \dz \bigg)^{\frac{p(n+1)-n}{p(n+2)}}\\
    &=
    \frac1M 
    \bigg(\int_{2^5Q_i^+} (M|f|)^{p_\#'} \boldsymbol{F}\dz \bigg)^{\frac{p(n+1)-n}{p(n+2)}}\\
    &\le
    \frac1M\,(\lambda^p |Q_i^+|)^{\frac{p(n+1)-n}{p(n+2)}} .
\end{align*}
Inserting the estimates for IV$_1$ and IV$_2$ above we get
\begin{align*}
    \mbox{IV}
    \le
    c(n,N,p,\nu,L,\psi,\|\omega\|_\infty)\,
    \frac{\lambda^p}{M}\, |Q_i^+| .
\end{align*}
Joining the preceding estimates for the terms I -- IV and \eqref{compare-r2} with \eqref{compare-1} we obtain the following {\it comparison estimate} for $Dv_i$:
\begin{equation}\label{compare-v}
    \int_{2^5Q_i^+} 
    \big(\mu^2 + |Dv_i|^2 + |Du|^2\big)^{\frac{p-2}{2}} |Dv_i-Du|^2 \dz \le
    c \bigg[\delta  +
    \frac{1}{M} +
    \frac{1}{M^{p-1}} \bigg]
    \lambda^p\,|Q_i^+| 
\end{equation}
with a constant $c=c(n,N,p,\nu,L,\psi,\|\omega\|_\infty)$.

In the following, we consider a second comparison problem in order to ``freeze'' the vector field $b$ with respect to the spatial variable $x$ and in the case that $b$ depends on $u$ also with respect to $u$. Here, we abbreviate
\begin{equation*}
    \mathcal B_i(t,\xi)
    :=
    \big(b(\cdot,t,(u+g)_{2^4Q_i^+},\xi)\big)_{2^4 B_i^+}
    \quad\mbox{for $t\in2^4\Lambda_i$ and $\xi\in\R^{Nn}$.}
\end{equation*}
By
\begin{equation*}
    w_i\in
    C^0\big(2^4\Lambda_i;L^2(2^4 B_i^+,\R^N)\big)
    \cap
    L^p\big(2^4\Lambda_i;W^{1,p}(2^4 B_i^+,\R^N)\big)
\end{equation*}
we denote the unique solution the to following Cauchy-Dirichlet problem
\begin{equation}\label{system-compare-2}
	\left\{
	\begin{array}{cc}
	\partial_t w_i -
    \Div \mathcal B_i(t,Dw_i)
	=
	0
	\quad &
	\text{in } 2^4 Q_i^+\\[7pt]
	w_i=v_i
	\quad &
	\text{on } \partial_p 2^4 Q_i^+.
	\end{array}
	\right.
\end{equation}
To infer energy and comparison estimates for $w_i$ we proceed similarly as before.
Subtracting the weak forms of \eqref{system-compare-2} and \eqref{system-compare-1} we get
\begin{align*}
    -\int_{2^4 Q_i^+}
    (w_i-v_i)\cdot\partial_t\varphi\dz+
    \int_{2^4 Q_i^+} & \big\langle \mathcal B_i(t,Dw_i),
    D\varphi\big\rangle \dz \nn\\
    & =
    \int_{2^4 Q_i^+}
    b(z,u+g,Dv_i\,\Psi), D\varphi\,\Psi\big\rangle \dz
\end{align*}
whenever $\varphi\in C_0^\infty(2^4 Q_i^+, \R^N)$.
For $\tau\in2^4\Lambda_i$ and $\epsilon>0$ such that $\tau+\epsilon\in2^4\Lambda_i$ we now choose the testing-function $\varphi_\epsilon=(w_i-v_i)\chi_\epsilon$, where $\chi_\epsilon$ is defined in \eqref{chi-epsilon}. Then, as in \eqref{time-part} we infer that the first integral on the right-hand side is non-negative. Hence, passing to the limit $\epsilon\downarrow 0$ and letting $\tau\uparrow t_{i}+\lambda^{2-p}(2^4\rho_{\mathfrak z_i})^2$ we obtain
\begin{align}\label{test-diff-w}
    \int_{2^4 Q_i^+} & \big\langle \mathcal B_i(t,Dw_i),
    D(w_i-v_i)\big\rangle \dz \\
    & \le
    \int_{2^4 Q_i^+} \big\langle b(z,u+g,Dv_i\,\Psi),D(w_i-v_i)\Psi\big\rangle \dz.
    \nn
\end{align}
This inequality will be used in two directions. First, we shall derive an energy bound for $Dw_i$.
Rearranging terms in \eqref{test-diff-w} and using \eqref{monotone-2}, \eqref{assump-b}$_1$ and Young's inequality we obtain from \eqref{test-diff-w} that for $\vartheta>0$ there holds
\begin{align*}
    \frac{\nu}{c(n,N,p)} \int_{2^4 Q_i^+} |Dw_i|^p \dz
    &\le
    \int_{2^4 Q_i^+}\big\langle \mathcal B_i(t,Dw_i),
    Dw_i\big\rangle + \mu^p \dz \nn\\
    & \le
    L\int_{2^4 Q_i^+} \big(\mu^2 + |Dw_i|^2\big)^{\frac{p-1}{2}} |Dv_i| + \mu^p\dz \\
    &\phantom{\le\ }+
    L \int_{2^4 Q_i^+} \big(\mu^2 + |Dv_i\,\Psi|^2\big)^{\frac{p-1}{2}} |D(w_i-v_i)\Psi| \dz \nn\\
    & \le
    \vartheta\int_{2^4 Q_i^+} |Dw_i|^p \dz +
    c\int_{2^4 Q_i^+} \big(|Dv_i|^p + 1\big) \dz,
\end{align*}
where $c=c(n,N,p,L,\psi,1/\vartheta)$.
Choosing $\vartheta$ small enough, i.e. $\vartheta=\nu/(2c(n,N,p))$ we can re-absorb the first integral of the right-hand side into the left. Subsequently using \eqref{energy-v} we obtain the following {\it energy estimate for $Dw_i$}:
\begin{align}\label{energy-w}
    \int_{2^4 Q_i^+} |Dw_i|^p \dz
    &\le
    c \int_{2^4 Q_i^+} \big(|Dv_i|^p + 1\big) \dz \\
    &\le
    c(n,N,p,\nu,L,\psi,\|\omega\|_\infty)\,\lambda^p\,|Q_i^+|.
    \nn
\end{align}
We now come to the comparison estimate for $w_i$.
Here, we again start from \eqref{test-diff-w}. Rearranging terms we get
\begin{align}\label{compare-2}
    \int_{2^4 Q_i^+} & \big\langle\mathcal B_i(t,Dw_i) - \mathcal B_i(t,Dv_i),
    D(w_i-v)\big\rangle \dz \\
    &\le
    \int_{2^4 Q_i^+}
    \big\langle b(z,u+g,Dv_i\,\Psi),
    D(w_i-v_i)(\Psi - \mathbb I_{n\times n})\big\rangle \dz \nn\\
    &\phantom{\le\ }+
    \int_{2^4 Q_i^+}
    \big\langle b(z,u+g,Dv_i\,\Psi) - b\big(z,u+g,Dv_i),
    D(w_i-v_i)\big\rangle \dz \nn\\
    &\phantom{\le\ }+
    \int_{2^4 Q_i^+}
    \big\langle b(z,u+g,Dv_i) - b(z,(u+g)_{2^4 Q_i^+},Dv_i), D(w_i-v_i)\big\rangle \dz \nn\\
    &\phantom{\le\ }+
    \int_{2^4 Q_i^+}
    \big\langle b(z,(u+g)_{2^4 Q_i^+},Dv_i) - \mathcal B_i(t,Dv_i) ,D(w_i-v_i)\big\rangle \dz \nn\\
    &=:
    \mbox{I} + \mbox{II} + \mbox{III} + \mbox{IV}.
    \nn
\end{align}
For the integral on the left-hand side of \eqref{compare-2} we get by \eqref{monotone} the following lower bound:
\begin{align*}
    \int_{2^4 Q_i^+} & \big\langle\mathcal B_i(t,Dw_i) - \mathcal B_i(t,Dv_i),
    D(w_i-v_i)\big\rangle \dz \nn\\
    &\ge
    \frac{\nu}{c(n,N,p)} \int_{2^4 Q_i^+}
    \big(\mu^2 + |Dw_i|^2 + |Dv_i|^2\big)^{\frac{p-2}{2}} |Dw_i-Dv_i|^2 \dz.
\end{align*}
From \eqref{assump-b}$_1$, the facts that $0\in Q_{2R}^+(z_o)$ (see \eqref{system-lat-CZ}), $2^4 Q_i^+\subset Q_{2R}^+(z_o)$ and $\Psi(0)=\mathbb I_{n\times n}$ by \eqref{Psi-CZ}, Young's inequality and the energy estimates \eqref{energy-v} and \eqref{energy-w} we get
\begin{align*}
    \mbox{I}
    &\le
    \int_{2^4 Q_i^+}
    |b(z,u+g,Dv_i\,\Psi)|\, |D(w_i-v_i)| |\Psi(x)-\mathbb I_{n\times n}| \dz \nn\\
    &\le
    L \osc_{B_{2R}^+(x_o)}(\Psi)
    \int_{2^4 Q_i^+} \big(\mu^2 + |Dv_i\,\Psi|^2\big)^{\frac{p-1}{2}} |Dw_i - Dv_i| \dz \nn\\
    &\le
    c(p)\,L \psi^{p-1}\osc_{B_{2R}^+(x_o)}(\Psi)
    \int_{2^4 Q_i^+} \big(\mu^p + |Dv_i|^p + |Dw_i|^p\big) \dz \\
    &\le
    c(n,N,p,\nu,L,\psi,\|\omega\|_\infty)\, \osc_{B_{2R}^+(x_o)}(\Psi)\,
    \lambda^p\,|Q_i^+|.
\end{align*}
For the estimate of II we use \eqref{assump-b}$_1$, Lemma \ref{lem:Fusco}, $0\in Q_{2R}^+(z_o)$, $\Psi(0)=\mathbb I_{n\times n}$ and the energy estimates \eqref{energy-v} and \eqref{energy-w} to infer that
\begin{align*}
    \mbox{II}
    &=
    \int_{2^4 Q_i^+}\int_0^1 
    \!\big\langle D_\xi b\big(z,u{+}g,Dv_i + sDv_i(\Psi {-} \mathbb I_{n\times n})\big)
    Dv_i(\Psi {-} \mathbb I_{n\times n}),
    D(w_i-v_i)\big\rangle \ds\dz \\
    &\le
    L\int_{2^4 Q_i^+}\int_0^1 
    \!\big(\mu^2 + |Dv_i {+} sDv_i(\Psi {-} \mathbb I_{n\times n})|^2\big)^{\frac{p-2}{2}} ds\,
    |Dv_i(\Psi {-} \mathbb I_{n\times n})|
    |D(w_i-v_i)| \dz \\
    &\le
    c(p)\,L \int_{2^4 Q_i^+}
    \!\big(\mu^2+|Dv_i|^2 + |Dv_i(\Psi{-}\mathbb I_{n\times n})|^2\big)^{\frac{p-2}{2}}
    |Dv_i(\Psi{-}\mathbb I_{n\times n})| 
    |D(w_i-v_i)| \dz \\
    &\le
    c(p)\,L\,\psi^{p-1} \osc_{B_{2R}^+(x_o)}(\Psi)
    \int_{2^4 Q_i^+} \big(\mu^2+|Dv_i|^2\big)^{\frac{p-2}{2}} |Dv_i| |Dw_i-Dv_i| \dz \\
    &\le
    c(p)\, \osc_{B_{2R}^+(x_o)}(\Psi)
    \int_{2^4 Q_i^+} \big(\mu^p + |Dv_i|^p + |Dw_i|^p\big) \dz \\
    &\le
    c(n,N,p,\nu,L,\psi,\|\omega\|_\infty)\, \osc_{B_{2R}^+(x_o)}(\Psi) \,
    \lambda^p\,|Q_i^+| .
\end{align*}
We now come to the estimate for the third term of \eqref{compare-2}.
Here, we first note that $\mbox{III} = 0$ if the vector field $b$ is independent of $u$, i.e. if $b(z,u,\xi)=b(z,\xi)$, while in the case where $b(z,u,\xi)$ depends on $u$ we assume that $u,g\in C^{0}(Q^+_{2R}(z_o)\cup\Gamma_{2R}(z_o),\R^N)$.
In the latter case we get with the help of assumption \eqref{assump-b}$_3$, the energy estimates \eqref{energy-v} and \eqref{energy-w} that
\begin{align*}
    \mbox{III}
    &\le
    \int_{2^4 Q_i^+}
    \theta\big(|u+g - (u+g)_{2^4 Q_i^+}|\big) \big(\mu^2+|Dv_i|^2\big)^{\frac{p-1}{2}} |D(w_i-v_i)| \dz \\
    &\le
    \psi^{p-1}\,
    \theta\big(\osc_{Q_{2R}^+(x_o)}(u) + \osc_{Q_{2R}^+(x_o)}(g)\big)
    \int_{2^4 Q_i^+} \big(\mu^2 + |Dv_i|^2\big)^{\frac{p-1}{2}} |Dw_i-Dv_i| \dz\\
    &\le
    c\, \theta\big(\osc_{Q_{2R}^+(x_o)}(u) + \osc_{Q_{2R}^+(x_o)}(g)\big)
    \int_{2^4 Q_i^+} \big(\mu^p + |Dv_i|^p + |Dw_i|^p\big) \dz \\
    &\le
    c(n,N,p,\nu,L,\psi,\|\omega\|_\infty)\, 
    \theta\big(\osc_{Q_{2R}^+(x_o)}(u) + \osc_{Q_{2R}^+(x_o)}(g)\big) \,
    \lambda^p\,|Q_i^+| .
\end{align*}
For the estimate of IV we use the definition of $\mathcal B_i$ and the VMO-assumption \eqref{VMO-b} which yield that
\begin{align*}
    \mbox{IV}
    &\le
    \int_{2^4 Q_i^+} \mint_{2^4 B_i^+} 
    \big|b\big(x,t,(u+g)_{2^4 Q_i^+},Dv_i\big) \\
    &\phantom{mmmmmmmmmmm} - 
    b\big(y,t,(u+g)_{2^4 Q_i^+},Dv_i\big)\big| \dy\, 
    |Dw_i-Dv_i| \dz \\
    &\le
    \int_{2^4 Q_i^+}
    {\bf v}_{\mathfrak x_i}(x,2^4\rho_{\mathfrak z_i}) (1+|Dv_i|)^{p-1} |Dw_i-Dv_i| \dz \\
    &\le
    \bigg(\int_{2^4 Q_i^+} {\bf v}_{\mathfrak x_i}(x,2^4\rho_{\mathfrak z_i})^{p'} (1+|Dv_i|)^{p} \dz \bigg)^{\frac1{p'}}
    \bigg(\int_{2^4 Q_i^+} |Dw_i-Dv_i|^p \dz \bigg)^{\frac1p}
\end{align*}
provided 
\begin{equation}\label{vmo-radius}
 	2^4\rho_{\mathfrak z_i}\le R\le R_o\le\rho_o.
\end{equation}
To further estimate the first integral on the right-hand side of the preceding inequality we use H\"older's inequality, the fact that ${\bf v}_{x_i}\le 2L$ and the higher integrability result from Lemma \ref{lem:hi} -- which is applicable due to \eqref{energy-v} -- to infer that
\begin{align*}
    & \int_{2^4 Q_i^+} {\bf v}_{\mathfrak x_i}(x,\rho/2)^{p'} (1+|Dv_i|)^{p} \dz  \\
    &\phantom{n}\le
    |2^4 Q_i^+|\bigg(\mint_{2^4 B_i^+} {\bf v}_{\mathfrak x_i}(x,2^4\rho_{\mathfrak z_i})^{\frac{p'(1+\epsilon_o)}{\epsilon_o}} \dx\bigg)^{\frac{\epsilon_o}{1+\epsilon_o}}
    \bigg(\mint_{2^4 Q_i^+} (1+|Dv_i|)^{p(1+\epsilon_o)} \dz \bigg)^{\frac1{1+\epsilon_o}} \\
    &\phantom{n}\le
    c \bigg(\mint_{2^4 B_i^+} {\bf v}_{\mathfrak x_i}(x,2^4\rho_{\mathfrak z_i})^{\frac{2p'}{\epsilon_o}} \dx\bigg)^{\frac{\epsilon_o}2}
    \lambda^p\,|Q_i^+| \\
    &\phantom{n}\le
    c \bigg(\mint_{2^4 B_i^+} {\bf v}_{\mathfrak x_i}(x,2^4\rho_{\mathfrak z_i}) \dx\bigg)^{\frac{\epsilon_o}{2}}
    \lambda^p\,|Q_i^+| \\
    &\phantom{n}\le
    c(n,N,p,\nu,L,\psi, \|\omega\|_\infty)\,\V(2R)^{\frac{\epsilon_o}{2}}\,
    \lambda^p\, |Q_i^+| .
\end{align*}
For the second integral we have by the energy estimates \eqref{energy-w} and \eqref{energy-v} that
\begin{align*}
    \int_{2^4 Q_i^+} |Dw_i-Dv_i|^p \dz
    &\le
    c(p) \int_{2^4 Q_i^+} \big(|Dw_i|^p + |Dv_i|^p\big) \dz \\
    & \le
    c(n,N,p,\nu,L,\psi, \|\omega\|_\infty)\,
    \lambda^p\, |Q_i^+|.
\end{align*}
Inserting the preceding two estimates above we get
\begin{align*}
    \mbox{IV}
    &\le
    c(n,N,p,\nu,L,\psi, \|\omega\|_\infty)\,\V(2R)^{\frac{\epsilon_o}{2 p'}}\,
    \lambda^p\, |Q_i^+| .
\end{align*}
Joining the previous estimates for I -- IV with \eqref{compare-2} we finally obtain the {\it comparison estimate} for $Dw_i$:
\begin{align}\label{compare-w}
    \int_{2^4 Q_i^+} & 
    \big(\mu^2+ |Dw_i|^2+ |Dv_i|^2\big)^{\frac{p-2}{2}}  |Dw_i-Dv_i|^2 \dz \\
    &\phantom{mmmmmmmm}\le
    c\Big[\widetilde\omega(2R) + \V(2R)^{\frac{\epsilon_o}{2p'}}\Big]\,
    \lambda^p\, |Q_i^+| \nn
\end{align}
with a constant $c$ depending on $n,N,p,\nu,L,\psi, \|\omega\|_\infty$. 
Here, we have abbreviated 
\begin{equation*}
	\widetilde\omega(2R)
	:=
	\osc_{B_{2R}^+(x_o)}(\Psi) +\theta\big(\osc_{Q_{2R}^+(x_o)}(u) + \osc_{Q_{2R}^+(x_o)}(g)\big).
\end{equation*}
At this point we recall that $\Psi$ is continuous on $B_{2R}^+(x_o)$ and that
$u$ and $g$ are assumed to be continuous on $Q_{2R}^+(z_o)$ when the vector field $b$ depends on $u$. Otherwise we may assume $\theta(\cdot)\equiv 0$. In any case we have
\begin{equation}\label{tilde-omega}
	\lim_{r\downarrow 0}\widetilde\omega(r)
	=
	0.
\end{equation}
Finally, we combine \eqref{compare-v} and \eqref{compare-w} with the help of Lemma \ref{lem:V} (ii). This yields 
\begin{align}\label{compare}
    \int_{2^4 Q_i^+} &
    \big(\mu^2+ |Dw_i|^2 + |Du|^2\big)^{\frac{p-2}{2}} |Dw_i-Du|^2 \dz \\
    &\le
    c\,\bigg[\delta +
    \frac{1}{M} +
    \frac{1}{M^{p-1}} +
    \widetilde\omega(2R) +
    \V(2R)^{\frac{\epsilon_o}{2p'}}\bigg]\,
    \lambda^p\, |Q_i^+|, \nn
\end{align}
where the constant $c$ depends on $n,N,p,\nu,L,\psi,\|\omega\|_\infty$.

\subsection{Estimates on intrinsic cylinders}
By (\ref{energy-w}) there exists a constant $c_\ast=c_\ast(n,N,p,$ $\nu,L,\psi,\|\omega\|_\infty)$ such that
\begin{align*}
    \mint_{2^4 Q_i^+} |Dw_i|^p \dz
    \le
    c_\ast\lambda^p.
\end{align*}
Therefore, we can use assumption \eqref{w-cond} to infer the existence of $H_b\equiv H_b(c_\ast,\psi, b(\cdot))\ge 1$ such that
\begin{align*}
    \mint_{2^3 Q_i^+} |Dw_i|^\chi \dz
    \le
    H_b\,\lambda^\chi.
\end{align*}
Recalling that $\lambda\ge 1$ and $\mu\in[0,1]$ we thus obtain
\begin{align}\label{aux-est}
    \mint_{2^3 Q_i^+} \big(\mu^2 + |Dw_i|^2\big)^{\frac{\chi}{2}} \dz
    \le
    2^\chi H_b\,\lambda^\chi.
\end{align}
For $A \ge 1$ to be chosen later and $z\in 2^3Q_i^+$ such that 
$|D u(z)| > A\lambda$ we will now derive a suitable estimate for $|Du(z)|^p$.
We first recall from Lemma \ref{lem:V} (i) that
\begin{align*}
	|D u(z)|^p
	&\le
	c_\ell\,\big(\mu^2 + |D w_i(z)|^2\big)^{\frac{p}{2}} +
	c_\ell\,\big(\mu^2 + |D w_i(z)|^2 + |D u(z)|^2\big)^{\frac{p-2}{2}}
	|D w_i-D u|^2  ,
\end{align*}
where $c_\ell \equiv c_\ell(n,N,p)$ is the constant from Lemma \ref{lem:V}.
We now distinguish two cases. In the first case where
\begin{align}\label{cases}
	\big(\mu^2 + |D w_i(z)|^2\big)^{\frac{p}{2}} 
	\le 
	\big(\mu^2 + |D w_i(z)|^2 + |D u(z)|^2\big)^{\frac{p-2}{2}} |D w_i(z)-D u(z)|^2  
\end{align}
we have 
\begin{align*}
	|D u(z)|^p
	&\le
	2c_\ell\,\big(\mu^2 + |D w_i(z)|^2 + |D u(z)|^2\big)^{\frac{p-2}{2}}
	|D w_i(z)-D u(z)|^2  ,
\end{align*}
while in the second case where the contrary of \eqref{cases} is true there holds
\begin{align*}
	A^p\lambda^p
	<
	|D u(z)|^p
	&\le
	2 c_\ell\,\big(\mu^2 + |D w_i(z)|^2\big)^{\frac{p}{2}} 
\end{align*}
and hence 
\begin{align*}
	|D u(z)|^p
	&\le
	2 c_\ell\,\big(\mu^2 + |D w_i(z)|^2\big)^{\frac{p}{2}} \cdot
	\Big[\frac{2 c_\ell}{A^p\lambda^p}\,\big(\mu^2 + |D w_i(z)|^2\big)^{\frac{p}{2}}\Big]^{\frac{\chi}{p}-1} \\
	&=
	\frac{(2 c_\ell)^{\frac\chi p}}{(A\lambda)^{\chi-p}} 
	\,\big(\mu^2 + |D w_i(z)|^2\big)^{\frac{\chi}{2}} .
\end{align*}
Therefore, in any case we have 
\begin{align*}
	|D u(z)|^p
	&\le
	2c_\ell\,\big(\mu^2 + |D w_i(z)|^2 + |D u(z)|^2\big)^{\frac{p-2}{2}}
	|D w_i(z)-D u(z)|^2  \\
	&\phantom{\le\ }+
	\frac{(2 c_\ell)^{\frac\chi p}}{(A\lambda)^{\chi-p}} 
	\,\big(\mu^2 + |D w_i(z)|^2\big)^{\frac{\chi}{2}} ,
\end{align*}
provided $|D u(z)| > A\lambda$. This implies the following estimate
\begin{align*}
	\int_{2^3Q_i^+ \cap E(A\lambda,r_2)} | D u|^p \dz
	&\le
	2c_\ell 
	\int_{2^3Q_i^+} \big(\mu^2 + |D w_i|^2 + |Du|^2\big)^{\frac{p-2}{2}}
	|D w_i-D u|^2 \dz  \\
	&\phantom{\le\ }+
	\frac{(2 c_\ell)^{\frac\chi p}}{(A\lambda)^{\chi-p}} 
	\int_{2^3Q_i^+} \big(\mu^2 + |D w_i|^2\big)^{\frac{\chi}{2}} \dz.
\end{align*}
Due to the comparison estimate \eqref{compare} we can bound the first term of the right-hand side of the preceding inequality while for the second one we use the gradient higher integrability bound \eqref{aux-est}. This leads us to
\begin{align*}
	\int_{2^3Q_i^+ \cap E(A\lambda,r_2)} | D u|^p \dz 
	&\le
	c\,\lambda^p \bigg[\delta +
    \frac{1}{M} + \frac{1}{M^{p-1}} +
    \widetilde\omega(2R) + \V(2R)^{\frac{\epsilon_o}{2p'}} 
    \bigg]
    |Q_i^+|  \\
	&\phantom{\le\ }+
	\frac{(2c_\ell)^{\frac{\chi}{p}}}{A^{\chi-p}} \, 2^{3(n+2)}2^\chi \lambda^p H_b\,
	|Q_i^+| .
\end{align*}
Note that $A,M\ge 1$ are yet not fixed and that $c$ depends on $n,N,p,\nu,L,\psi,\|\omega\|_\infty$.
Hence, setting
\begin{equation}\label{def-M}
	\mathcal M \equiv \mathcal M(\delta ,R,M,A)
	:=
	\delta +
    \frac{1}{M} + \frac{1}{M^{p-1}} +
    \widetilde\omega(2R) + \V(2R)^{\frac{\epsilon_o}{2p'}} +
    \frac{1}{A^{\chi-p}}
\end{equation}
we have
\begin{align}\label{Qi-level}
	\int_{2^3Q_i^+ \cap E(A\lambda,r_2)} | D u|^p \dz
	\le
	c\,\mathcal M\,\lambda^p |Q_i^+|
\end{align}
with a constant $c$ depending on $n,N,p,\nu,L,\psi,\|\omega\|_\infty, H_b$.

Next, we shall derive a suitable bound for $|Q_i^+|$. We first rewrite (\ref{stop-Qj}) for $j=0$ in the form
\begin{equation}\label{Q0-1}
	|Q_i^+|
	=
	\frac{1}{\lambda^p} \int_{Q_i^+} |D u|^p \dz +
	\frac{M^p }{\lambda^p} \int_{Q_i^+} \G^p \dz +
	\frac{M^{p_\#' } \boldsymbol{F} }{\lambda^p} \int_{Q_i^+} |f|^{p_\#'} \dz  .
\end{equation}
To estimate the first integral in (\ref{Q0-1}) we decompose the domain of integration into the set where
$|D u|>\lambda/8$, respectively $|D u|\le\lambda/8$,
to obtain
\begin{align*}
	\frac{1}{\lambda^p} \int_{Q_i^+} |D u|^p \dz
	%&=
%	\frac{1}{\lambda^p}
%	\int_{Q_i^{(0)}\cap\{|D u|>\tau\lambda\}} |D u|^p \dz +
%	\frac{1}{\lambda^p}
%	\int_{Q_i^{(0)}\cap\{|D u|\le\tau\lambda\}}
%	|D u|^p \dz \\
	&\le
	\frac{1}{\lambda^p}
	\int_{Q_i^+\cap\{|D u|>\lambda/8\}} |D u|^p \dz +
	8^{-p} |Q_i^+|.
\end{align*}
Similarly, the second integral in (\ref{Q0-1}) can be estimated by
\begin{align*}
	\frac{M^p}{\lambda^p} \int_{Q_i^+} \G^p \dz
	%&=
%	\frac{M^p}{\lambda^p}
%	\int_{Q_i^{(0)}\cap\{\Psi_\mu>\gamma\lambda\}} \Psi_\mu^p \dz +
%	\frac{M^p }{\lambda^p}
%	\int_{Q_i^{(0)}\cap\{\Psi_\mu\le\gamma\lambda\}}
%	\Psi^p \dz \\
	&\le
	\frac{M^p }{\lambda^p}
	\int_{Q_i^+\cap\{\G>\lambda/(8M)\}} \G^p \dz +
	8^{-p} |Q_i^+|
\end{align*}
and the third one by 
\begin{align*}
	\frac{M^{p_\#' } \boldsymbol{F} }{\lambda^p} \int_{Q_i^+} |f|^{p_\#'} \dz
	&\le
	\frac{M^{p_\#' } \boldsymbol{F} }{\lambda^p}
	\int_{Q_i^+\cap\{|f|^{p_\#'}>\lambda^p/(8M^{p_\#'}\boldsymbol{F})\}} |f|^{p_\#'} \dz +
	8^{-1} |Q_i^+|.
\end{align*}
Combining the preceding estimates with (\ref{Q0-1})
and reabsorbing $\frac{1}{2}|Q_i^+|$ on the left-hand side we obtain
the {\it measure estimate}
\begin{align}\label{Q0-2}
	|Q_i^+|
	&\le
	\frac{2}{\lambda^p} \,\bigg[
	\int_{Q_i^+\cap\{|D u|>\lambda/8\}} |D u|^p \dz +
	M^p\int_{Q_i^+\cap\{\G>\lambda/(8M)\}} \G^p \dz \\
	&\phantom{mmmmmmmmmm}+
	M^{p_\#' } \boldsymbol{F}
	\int_{Q_i^+\cap\{|f|^{p_\#'}>\lambda^p/(8M^{p_\#'}\boldsymbol{F})\}} |f|^{p_\#'} \dz\bigg].
	\nn
\end{align}
Connecting (\ref{Q0-2}) with (\ref{Qi-level}) yields
\begin{align}\label{Q0-3}
	& \int_{2^3Q_i^+ \cap E(A\lambda,r_2)} | D u|^p \dz \\
	&\phantom{mm}\le
	c\,\mathcal M
		\bigg[
    	\int_{Q_i^+\cap\{|D u|>\lambda/8\}} |D u|^p \dz +
	M^p \int_{Q_i^+\cap\{\G>\lambda/(8M)\}} \G^p\dz \nn\\
	&\phantom{mmmmmmmmmmmm}+
	M^{p_\#' } \boldsymbol{F}
	\int_{Q_i^+\cap\{|f|^{p_\#'}>\lambda^p/(8M^{p_\#'}\boldsymbol{F})\}} |f|^{p_\#'} \dz\bigg]
	\nonumber
\end{align}
with a constant $c$ depending on $n,N,p,\nu,L,\psi,\|\omega\|_\infty, H_b$.

\subsection{Concluding the proof}
We now derive an upper bound for the super-level sets $E(A\lambda, r_1)$ introduced in \eqref{def-levelset}.
Recalling that on the one hand the family $(2^3Q_i^+)_{i\in\N}$ covers $E(\lambda, r_1)\supset E (A\lambda,r_1)$ and on the other hand the corresponding family $(Q_i^+)_{i\in\N}$ consists of pairwise disjoint cylinders $Q_i^+\subset Q_{r_2}^+(z_o)$
we obtain after summing (\ref{Q0-3}) over $i\in\N$ that
\begin{align*}
	& \int_{E(A\lambda,r_1)} | D u|^p \dz \\
	&\phantom{mm}\le
	c\, \mathcal M
	\bigg[\int_{Q_{r_2}^+(z_o)\cap\{|D u|>\lambda/8\}} |D u|^p \dz +
	M^p\int_{Q_{r_2}^+(z_o)\cap\{\G>\lambda/(8M)\}} \G^p \dz \\
	&\phantom{mmmmmmmmmmmmmn}+
	M^{p_\#' } \boldsymbol{F}
	\int_{Q_{r_2}^+(z_o)\cap\{|f|^{p_\#'}>\lambda^p/(8M^{p_\#'}\boldsymbol{F})\}} |f|^{p_\#'} \dz\bigg].
\end{align*}
Note that $A,M\ge 1$ are still to be chosen and that by \eqref{def-B} and \eqref{lambda-asym} the preceding inequality is valid for any $\lambda\ge\lambda_1$, where
\begin{equation}\label{def-lambda1}
    \lambda_1
    :=
    \max\bigg\{ B\lambda_o\,,\, \frac{K_\delta}{\delta^{\frac{1}{p-1}}}\bigg\}\,.
\end{equation}
At this point we would like to multiply both sides of the previous inequality by $\lambda^{q-p-1}$ and then
integrate with respect to $\lambda$ over $[\lambda_1,\infty)$.
This, formally would lead to an $L^q$ estimate of $Du$ if certain parameters are chosen small enough to reabsorb $\int|Du|^q\dz$ on the left-hand side. However, we are not allowed to perform this step since the integral might be infinite. This problem will be overcome in the following by a truncation argument. For $k\ge \lambda_1$ we define
\begin{equation*}
    |Du|_k
    :=
    \min\{|Du|,k\}
\end{equation*}
and
\begin{equation*}
    E_k(A\lambda,r_1)
    :=
    \big\{z\in Q_{r_1}^+(z_o) :
    |Du(z)|_k >A\lambda\big\}.
\end{equation*}
Then, from the last inequality we deduce that
\begin{align*}
	& \int_{E_k(A\lambda,r_1)} | D u|^p \dz \\
	&\phantom{mm}\le
	c\, \mathcal M
	\bigg[\int_{Q_{r_2}^+(z_o)\cap\{|D u|_k>\lambda/8\}} |D u|^p \dz +
	M^p\int_{Q_{r_2}^+(z_o)\cap\{\G>\lambda/(8M)\}} \G^p \dz \\
	&\phantom{mmmmmmmmmmmmmm}+
	M^{p_\#' } \boldsymbol{F}
	\int_{Q_{r_2}^+(z_o)\cap\{|f|^{p_\#'}>\lambda^p/(8M^{p_\#'}\boldsymbol{F})\}} |f|^{p_\#'} \dz\bigg].
\end{align*}
This can be seen as follows: In the case $k\le A\lambda$ we have
$E_k(A\lambda,r_1 )=\emptyset$ and therefore the inequality holds trivially. In the second case $k>A\lambda$ the last inequality  follows from the second last one since then
$\{E_k(A\lambda,r_1 )\}= \{E(A\lambda,r_1 )\}$ and
$\{|D u|_k>\lambda/8\} = \{|D u|>\lambda/8\}$.
Now, we choose $q\in(p,\chi)$, multiply both sides of the preceding inequality by $\lambda^{q-1}$ and integrate with respect to $\lambda$ over $(\lambda_1,\infty)$. In this way we obtain
\begin{align}\label{level-1}
	\int_{\lambda_1}^\infty & 
	\lambda^{q-p-1} \int_{E_k(A\lambda,r_1)} |Du|^p \dz
	\,d\lambda \\
	&\le
	c\,\mathcal M
	\int_{\lambda_1}^\infty\lambda^{q-p-1} \bigg[
	\int_{Q_{r_2}^+(z_o)\cap\{|D u|_k>\lambda/8\}} |Du|^p
	\dz \nn\\
	&\phantom{mmmmmmmmmn}+
	 M^p\!\int_{Q_{r_2}^+(z_o)\cap\{\G>\lambda/(8M)\}}
	\G^p\dz \nn\\
	&\phantom{mmmmmmmmmn}+
	M^{p_\#' } \boldsymbol{F}
	\int_{Q_{r_2}^+(z_o)\cap\{|f|^{p_\#'}>\lambda^p/(8M^{p_\#'}\boldsymbol{F})\}} |f|^{p_\#'} \dz\bigg]
	d\lambda .
	\nonumber
\end{align}
To the integral on the left-hand side of \eqref{level-1} we apply Fubini's theorem which yields that
\begin{align*}
	\int_{\lambda_1}^\infty & 
	\lambda^{q-p-1} \int_{E_k(A\lambda,r_1)} |Du|^p \dz
	\,d\lambda \\
    &=
    \int_{E_k(A\lambda_1,r_1)} |Du|^p
    \int_{\lambda_1}^{|Du(z)|_k /A}
    \lambda^{q-p-1} \,d\lambda \dz
    \nonumber\\
    &=
    \frac{1}{q-p} \int_{E_k(A\lambda_1,r_1)} |Du|^p \Big[
    A^{-(q-p)}|Du|_k^{q-p} - \lambda_1^{q-p} \Big] \dz
    \nonumber\\
    &=
    \frac{1}{q-p} \bigg[
    \frac{1}{A^{q-p}} \int_{E_k(A\lambda_1,r_1)}
    |Du|^{p} |Du|_k^{q-p} \dz -
    \lambda_1^{q-p} 
    \int_{E_k(A\lambda_1,r_1)} |Du|^{p} \dz \bigg].
\end{align*}
Next, we consider the terms on the right-hand side of \eqref{level-1}. For the first one we obtain, again by Fubini's theorem, that
\begin{align*}
	\int_{\lambda_1}^\infty & \lambda^{q-p-1} 
	\int_{Q_{r_2}^+(z_o)\cap\{|Du|_k>\lambda/8\}} |D u|^p
	\dz\,d\lambda\\
	&=
    \int_{Q_{r_2}^+(z_o)\cap\{|Du|_k>\lambda_1/8\}} |D u|^p
	\int_{\lambda_1}^{8|Du(z)|_k} \lambda^{q-p-1} \,d\lambda \dz\\
	&\le
	\frac{8^{q-p}}{q-p}
    \int_{Q_{r_2}^+(z_o)} |D u|^{p} |Du|_k^{q-p} \dz.
\end{align*}
Similarly, we get for the second term
\begin{align*}
	\int_{\lambda_1}^\infty & \lambda^{q-p-1} 
	\int_{Q_{r_2}^+(z_o)\cap\{\G>\lambda/(8M)\}} \G^p \dz
	\,d\lambda \\
	&=
	\int_{Q_{r_2}^+(z_o)\cap\{\G>\lambda_1/(8M)\}} \G^p
    \int_{\lambda_1}^{8M\G(z)}\lambda^{q-p-1}
	\,d\lambda \dz \\
	&\le
	\frac{(8M)^{q-p}}{q-p} \int_{Q_{r_2}^+(z_o)} \G^{q} \dz
\end{align*}
and also for the third term
\begin{align*}
	\int_{\lambda_1}^\infty & \lambda^{q-p-1} 
	\int_{Q_{r_2}^+(z_o)\cap\{|f|^{p_\#'}>\lambda^p/(8M^{p_\#'}\boldsymbol{F})\}} |f|^{p_\#'} \dz
	\,d\lambda \\
	&=
	\int_{Q_{r_2}^+(z_o)\cap\{|f|^{p_\#'}>\lambda_1^p/(8M^{p_\#'}\boldsymbol{F})\}} |f|^{p_\#'}
    	\int_{\lambda_1}^{(8\boldsymbol{F})^{1/p} (M|f(z)|)^{p_\#'/p}}\lambda^{q-p-1}
	\,d\lambda \dz \\
	&\le
	\frac{M^{\frac{p_\#'(q-p)}{p}} (8\boldsymbol{F})^{\frac{q}{p}-1}}{q-p} \int_{Q_{r_2}^+(z_o)} |f|^{\frac{p_\#' q}{p}} \dz.
\end{align*}
Connecting the preceding estimates with (\ref{level-1}) we get
\begin{align*}
	\int_{E_k(A\lambda_1,r_1)} &
    	|Du|_k^{q-p} |Du|^{p} \dz \nn\\
	&\le
	A^{q-p}\lambda_1^{q-p} 
    \int_{Q_{r_1}^+(z_o)} |Du|^{p} \dz \nn\\
    &\phantom{\le\ }+
	c\,A^{q-p}\, \mathcal M
	\bigg[\int_{Q_{r_2}^+(z_o)} |Du|_k^{q-p} |Du|^{p} \dz
	+
	M^{q} \int_{Q_{r_2}^+(z_o)} \G^{q} \dz \nn\\
	&\phantom{mmmmmmmmmmmmmmm}+
	M^{\frac{qp_\#'}{p}} \boldsymbol{F}^{\frac{q}{p}} \int_{Q_{r_2}^+(z_o)} |f|^{\frac{p_\#' q}{p}} \dz
	\bigg] .
\end{align*}
Taking into account that
\begin{align*}
	\int_{Q_{r_1}^+(z_o)\setminus E_k(A\lambda_1,r_1)} 
    	|Du|_k^{q-p} |Du|^{p} \dz 
	\le
	(A\lambda_1)^{q-p} \int_{Q_{r_1}^+(z_o)} |Du|^{p} \dz 
\end{align*}
we arrive at 
\begin{align}\label{CZ-1}
	\int_{Q_{r_1}^+(z_o)} &
    	|Du|_k^{q-p} |Du|^{p} \dz \\
	&\le
	A^{q-p}\lambda_1^{q-p} 
    \int_{Q_{r_1}^+(z_o)} |Du|^{p} \dz \nn\\
    &\phantom{\le\ }+
	c\,A^{q-p}\, \mathcal M
	\bigg[\int_{Q_{r_2}^+(z_o)} |Du|_k^{q-p} |Du|^{p} \dz
	+
	M^{q} \int_{Q_{r_2}^+(z_o)} \G^{q} \dz \nn\\
	&\phantom{mmmmmmmmmmmmmmm}+
	M^{\frac{qp_\#'}{p}} \boldsymbol{F}^{\frac{q}{p}} \int_{Q_{r_2}^+(z_o)} |f|^{\frac{p_\#' q}{p}} \dz
	\bigg] , \nn
\end{align}
where $\mathcal M\equiv \mathcal M(\delta ,R,M,A)$ was defined in \eqref{def-M} and the constant $c$ depends on $n,N,p,\nu,L,$ $\psi,\|\omega\|_\infty, H_b$.
We now perform the choices of parameters as follows.
We first choose $A\ge 1$ in dependence of $n,N,p,\nu,L,\chi,q,\psi, \|\omega\|_\infty,H_b$, large enough to have
\begin{equation*}
    c\,A^{-(\chi-q)}
    \le
    \frac18\, .
\end{equation*}
Next, we choose $M\ge 1$ in dependence of the same parameters large enough, to have
\begin{equation*}
	c\,A^{q-p}~
	\Big[\frac{1}{M} + \frac{1}{M^{(p-1)}}\Big]
	\le
	\frac{1}{8}
\end{equation*}
and $\delta>0$ still depending on the same parameters small enough to ensure
\begin{equation*}
	c\,A^{q-p}~\delta
	\le
	\frac{1}{8}\,.
\end{equation*}
Note that this also fixes $K_\delta$ in dependence of $\delta$ and $\omega$. 
Finally, we select $R_o\in(0,1]$ small enough to have
\begin{equation*}
	c\,A^{q-p}~
    \Big[ \widetilde\omega(2R_o) + \V(2R_o)^{\frac{\epsilon_o}{2p'}}\Big]
	\le
	\frac{1}{8}
	\quad\mbox{and}\quad
	R_o\le\rho_o .
\end{equation*}
Note that this is possible by \eqref{VMO-b} and \eqref{tilde-omega}. The latter condition ensures that \eqref{vmo-radius} is satisfied. With this choice, $R_o$ depends on $n,N,p,\nu,L,\chi,q,\psi,$ $ \|\omega\|_\infty, H_b, \V(\cdot)$ and the modulus of continuity of $\Psi$ and if the vector field $b$ depends on $u$ also on $\theta(\cdot)$ and the moduli of continuity of $u$ and $g$.
These choices allow us to bound the constant
in front of the integrals appearing on the right-hand side of (\ref{CZ-1}) by $\frac12$, i.e. we have
\begin{equation*}
	c\,A^{q-p}~\mathcal M
	\le
	\frac{1}{2},
\end{equation*}
so that (\ref{CZ-1}) turns into
\begin{align*}
	\int_{Q_{r_1}^+(z_o)} |Du|_k^{q-p} |Du|^{p} \dz
	&\le
	\frac{1}{2}
	\int_{Q_{r_2}^+(z_o)} |Du|_k^{q-p} |Du|^{p} \dz \\
	&+
	c\, \lambda_1^{q-p} 
    \int_{Q_{2R}^+(z_o)} |Du|^{p} \dz +
	c  \int_{Q_{2R}^+(z_o)} \G^{q} + \boldsymbol{F}^{\frac qp} |f|^{\frac{p_\#' q}{p}} \dz ,
\end{align*}
where $c$ depends on 
$n,N,p,\nu,L,\chi,q,\psi,\|\omega\|_\infty, H_b$.
Recalling the definition of $\lambda_1$ from \eqref{def-lambda1} and \eqref{def-B} and taking averages we obtain from the preceding estimate
\begin{align*}
	\mint_{Q_{r_1}^+(z_o)} |Du|_k^{q-p} |Du|^{p} \dz
	&\le
    	\frac{1}{2}
	\mint_{Q_{r_2}^+(z_o)} |Du|_k^{q-p} |Du|^{p} \dz \\
	& \phantom{\le\ }+
	c\, \Big(\frac{R}{r_2{-}r_1}\Big)^\beta
    	\lambda_o^{q-p} \mint_{Q_{2R}^+(z_o)} |Du|^{p} \dz \\
    	& \phantom{\le\ }+
    	c \mint_{Q_{2R}^+(z_o)}\G^{q} + \boldsymbol{F}^{\frac qp} |f|^{\frac{p_\#' q}{p}} \dz +c ,
\end{align*}
where $\beta :=(n+2)(q-p)d/p$ and $c$ now depends on 
$n,N,p,\nu,L,\chi,q,\psi,\omega(\cdot), H_b$.
Since $R\le r_1<r_2\le2R$ are arbitrary we can apply Lemma \ref{lem:it}
with the choices
\begin{equation*}
    	\phi( r)=\mint_{Q_{r}^+(z_o)} |Du|_k^{q-p} |Du|^{p} \dz,
    	\quad
    	\vartheta =\tfrac12,
    	\quad
    	A= c\,R^\beta\, \lambda_o^{q-p} \mint_{Q_{2R}^+(z_o)} |Du|^{p} \dz
\end{equation*}
and
\begin{equation*}
    	B
    	= 
    	c \mint_{Q_{2R}^+(z_o)}\G^{q} + \boldsymbol{F}^{\frac qp} |f|^{\frac{p_\#' q}{p}} \dz + c.
\end{equation*}
This leads us to
\begin{align*}
	\mint_{Q_{R}^+(z_o)} & |Du|_k^{q-p} |Du|^{p} \dz \\
	&\le
	c\, \lambda_o^{q-p} \mint_{Q_{2R}^+(z_o)} |Du|^{p} \dz +
    	c \mint_{Q_{2R}^+(z_o)}\G^{q} + \boldsymbol{F}^{\frac qp} |f|^{\frac{p_\#' q}{p}} \dz + c.
\end{align*}
Letting $k\to \infty$ which is possible by Fatou's lemma we get
\begin{align*}
	\mint_{Q_{R}^+(z_o)} |Du|^q \dz
	&\le
	c\, \lambda_o^{q-p} \mint_{Q_{2R}^+(z_o)} |Du|^{p} \dz +
    	c \mint_{Q_{2R}^+(z_o)}\G^{q} + \boldsymbol{F}^{\frac qp} |f|^{\frac{p_\#' q}{p}} \dz + c.
\end{align*}
At this point we recall the definitions of $\lambda_o, \G,\boldsymbol{F}$ from \eqref{def-lambda0} and \eqref{def-G} and the fact that $d_{CZ}\ge 1$ (recall that $d_{CZ}$ is defined in \eqref{def-d}). This together with H\"older's inequality leads us to the final estimate
\begin{align*}
	\mint_{Q_{R}^+(z_o)} |Du|^q \dz
	&\le
	c \bigg[ \bigg( \mint_{Q_{2R}^+(z_o)} |D u|^p \dz \bigg)^\frac{q}{p}+
	\mint_{Q_{2R}^+(z_o)}\G^{q} + \boldsymbol{F}^{\frac qp} |f|^{\frac{p_\#' q}{p}} \dz + 1
	\bigg]^{d_{CZ}}  \\
	&\le
	c \bigg[ \bigg( \mint_{Q_{2R}^+(z_o)} |D u|^p \dz \bigg)^\frac{q}{p}+
	\mint_{Q_{2R}^+(z_o)} \big( |Dg|^{q} + |F|^q \big) \dz \\
	&\phantom{mmmm}+
	R^{-(n+2)} 
	\bigg(\int_{Q_{2R}^+(z_o)} |f|^{\frac{q(n+2)}{p(n+2)-n}} \dz \bigg)^{1+\frac{p}{np+p-n}} + 1
	\bigg]^{d_{CZ}} 
\end{align*}
with a constant $c$ depending on 
$n,N,p,\nu,L,\chi,q,\psi,\omega(\cdot), H_b$. 
This finishes the proof of Theorem \ref{thm:main-CZ}.
\hfill$\square$

\bibliographystyle{plain}

\begin{thebibliography}{10}

\bibitem{Acerbi-Fusco:1989}
E.~Acerbi and N.~Fusco,
\newblock Regularity for minimizers of nonquadratic functionals: the case $1<p<2$.
\newblock {\em J. Math. Anal. Appl.},  140(1):115--135, 1989.

\bibitem{Acerbi-Mingione:2005} 
E.~Acerbi and G.~Mingione,
\newblock Gradient estimates for the $p(x)$-Laplacean system. 
\newblock {\em J. Reine Angew. Math.}, 584:117--148, 2005.

\bibitem{Acerbi-Mingione:2007} 
E.~Acerbi and G.~Mingione,
\newblock Gradient estimates for a class of parabolic systems.
\newblock {\em Duke Math.~J.}, 136:285--320, 2007.

\bibitem{Boegelein:2008}
V.~B\"ogelein,
\newblock Higher integrability for weak solutions of higher order degenerate parabolic systems.
\newblock {\em Ann. Acad. Sci. Fenn., Math.}, 33(2):387--412, 2008.

\bibitem{Boegelein:2012}
V.~B\"ogelein,
\newblock The boundary regularity of nonlinear parabolic systems III.
\newblock {\em in preparation}.

\bibitem{Boegelein:2012-Lipschitz}
V.~B\"ogelein,
\newblock Global Lipschitz regularity for the parabolic $p$-Laplacian system.
\newblock {\em preprint}, 2012.

\bibitem{Boegelein-habil}
V.~B\"ogelein,
\newblock Global gradient estimates for degenerate and singular parabolic systems.
\newblock {\em Habilitationsschrift}, 2012.

\bibitem{Boegelein-Duzaar-Mingione:2010-boundary-II}
V.~B\"ogelein, F.~Duzaar and G.~Mingione,
\newblock The boundary regularity of nonlinear parabolic systems II.
\newblock {\em Ann. Inst. Henri Poincar\'e, Anal. Non Lin\'eaire}, 27:145--200, 2010.

\bibitem{Boegelein-Duzaar-Mingione:2011}
V.~B\"ogelein, F.~Duzaar and G.~Mingione,
\newblock Degenerate problems with irregular obstacles,
\newblock {\em J. Reine Angew. Math.}, 650:107--160, 2011.

\bibitem{Boegelein-Duzaar-Mingione:2010}
V.~B\"ogelein, F.~Duzaar and G.~Mingione,
\newblock The regularity of general parabolic systems with degenerate diffusions.
\newblock {\em Mem. Amer. Math. Soc.}, 221(1041).

\bibitem{Boegelein-Parviainen:2010}
V.~B\"ogelein and M.~Parviainen,
\newblock Self-improving property of nonlinear higher order parabolic systems near the boundary.
\newblock {\em Nonlinear Differ. Equ. Appl.}, 17(1):21--54, 2010.

%\bibitem{Byun:2007}
%S.~Byun and L.~Wang, 
%\newblock Quasilinear elliptic equations with BMO coefficients in Lipschitz domains. 
%\newblock {\em Trans. Amer. Math. Soc.}, 359(12):5899Ð5913, 2007.

\bibitem{Byun-Wang:2008}
S.~Byun and L.~Wang, 
\newblock Gradient estimates for elliptic systems in non-smooth domains. 
\newblock {\em Math. Ann.}, 341(3):629--650, 2008. 

\bibitem{Byun-Wang:2008-2}
S.~Byun and L.~Wang, 
\newblock Parabolic equations with BMO nonlinearity in Reifenberg domains. 
\newblock {\em J. Reine Angew. Math.}, 615:1--24, 2008.
35K55 (35K20)

\bibitem{Caffarelli-Peral:1998}
L.~Caffarelli and I.~Peral, 
\newblock On $W^{1,p}$ estimates for elliptic equations in divergence form. 
\newblock {\em Comm. Pure Appl. Math.}, 51:1--21, 1998.

\bibitem{Campanato:1984}
S.~Campanato,
\newblock On the nonlinear parabolic systems in divergence form. H\"older continuity and partial H\"older continuity of the solutions.
\newblock {\em Ann. Mat. Pura Appl. (IV)}, 137:83--122, 1984.

\bibitem{Chipot-Evans:1986}
M.~Chipot and L.~C.~Evans,
\newblock Linearisation at infinity and Lipschitz estimates for certain problems in the calculus of variations.
\newblock {\em Proc. Roy. Soc. Edinburgh Sect. A} 102(3-4):291--303, 1986.

\bibitem{Cupini-Fusco-Petti:1999}
G.~Cupini, N.~Fusco and R.~Petti,
\newblock H\"older continuity of local minimizers.
\newblock {\em J. Math. Anal. Appl.}, 235(2):578--597, 1999.

\bibitem{DiBenedetto:book} 
E.~DiBenedetto,
\newblock {\em Degenerate parabolic equations}.
\newblock Universitext. New York, NY: Springer-Verlag. xv 387, 1993.

\bibitem{DiBenedetto-Friedman:1984} 
E.~DiBenedetto and A.~Friedman,
\newblock Regularity of solutions of nonlinear degenerate parabolic systems. 
\newblock {\em J. Reine Angew. Math.}, 349:83--128, 1984.

\bibitem{DiBenedetto-Friedman:1985}
E.~DiBenedetto and A.~Friedman,
\newblock H\"older estimates for nonlinear degenerate parabolic systems.
\newblock{\em J. Reine Angew. Math.}, 357:1--22, 1985.

\bibitem{DiBenedetto-Manfredi:1993}
E.~DiBenedetto and J.~J.~Manfredi, 
\newblock On the higher integrability of the gradient of weak solutions of certain degenerate elliptic systems. 
\newblock {\em Amer. J. Math.}, 115:1107--1134, 1993.

\bibitem{Diening-Stroffolini-Verde:2011}
L.~Diening, B.~Stroffolini and A.~Verde, 
\newblock Lipschitz regularity for some asymptotically convex problems. 
\newblock {\em ESAIM Control Optim. Calc. Var.}, 17(1):178--189, 2011.

\bibitem{Duzaar-Mingione-Steffen:2011}
F.~Duzaar, G.~Mingione, and K.~Steffen,
\newblock Parabolic systems with polynomial growth and regularity,
\newblock {\em Mem. Amer. Math. Soc.}, 214(1005), 2011.

\bibitem{Foss:2008}
M.~Foss,
\newblock Global regularity for almost minimizers of nonconvex variational problems.
\newblock {\em Ann. Mat. Pura Appl. (4)}, 187(2):263--321, 2008.

\bibitem{Foss-Passarelli-Verde:2008}
M.~Foss, A.~Passarelli di Napoli and A.~Verde, 
\newblock Global Morrey regularity results for asymptotically convex variational problems. \newblock {\em Forum Math.}, 20(5):921--953, 2008.

\bibitem{Giaquinta-Modica:1986}
M.~Giaquinta and G.~Modica,
\newblock Remarks on the regularity of the minimizers of certain degenerate functionals.
\newblock {\em Manuscripta Math.}, 57(1):55--99, 1986.

\bibitem{Giusti:book}
E.~Giusti,
\newblock {\em Direct Methods in the Calculus of Variations}.
\newblock World Scientific Publishing Company, Tuck Link, Singapore, 2003.

\bibitem{Iwaniec:1983}
T.~Iwaniec, 
\newblock Projections onto gradient fields and $L^p$-estimates for degenerated elliptic operators.
\newblock {\em Studia Math.}, 75:293--312, 1983.

\bibitem{Kinnunen-Lewis:2000}
J.~Kinnunen and J.~L.~Lewis,
\newblock Higher integrability for parabolic systems of $p$-laplacian type.
\newblock {\em Duke Math. J.}, 102:253--271, 2000.

\bibitem{Kinnunen-Zhou-1999}
J.~Kinnunen and S.~Zhou,
\newblock A local estimate for nonlinear equations with discontinuous coefficients.
\newblock {\em Comm. Partial Differential Equations}, 24(11-12):2043--2068, 1999.

\bibitem{Kinnunen-Zhou-2001}
J.~Kinnunen and S.~Zhou,
\newblock A boundary estimate for nonlinear equations with discontinuous coefficients.
\newblock {\em Differential Integral Equations}, 14(4):475--492, 2001.

\bibitem{Kristensen-Mingione:2010}
J.~Kristensen and G.~Mingione, 
\newblock Boundary regularity in variational problems. 
\newblock {\em Arch. Ration. Mech. Anal.}, 198(2):369--455, 2010.

\bibitem{Misawa:2006}
M.~Misawa, 
\newblock $L^q$ estimates of gradients for evolutional $p$-Laplacian systems. 
\newblock {\em J. Differential Equations}, 219:390--420, 2006.

\bibitem{Parviainen:2009}
M.~Parviainen,
\newblock Global gradient estimates for degenerate parabolic equations in nonsmooth domains.
\newblock {\em Ann. Mat. Pura Appl.}, 188(2):333--358, 2009.

\bibitem{Raymond:1991}
J.~P.~Raymond,
\newblock Lipschitz regularity of solutions of some asymptotically convex problems.
\newblock {\em Proc. Roy. Soc. Edinburgh Sect. A}, 117(1-2):59--73, 1991.

\bibitem{Schemm:2009}
S.~Schemm,
\newblock Higher order functionals and regularity. 
\newblock Dissertation, 2009.

\bibitem{Scheven:2010}
C.~Scheven,
\newblock Regularity for subquadratic parabolic systems: higher integrability and dimension estimates. 
\newblock {\em Proc. Roy. Soc. Edinburgh Sect. A}, 140(6):1269--1308, 2010.

\bibitem{Scheven:2011}
C.~Scheven,
\newblock Existence of localizable solutions to nonlinear parabolic problems with irregular obstacles.
\newblock {\em preprint}, 2011.

\bibitem{Scheven-Schmidt:2009}
C.~Scheven and T. Schmidt,
\newblock Asymptotically regular problems. II. Partial Lipschitz continuity and a singular set of positive measure. 
\newblock {\em Ann. Sc. Norm. Super. Pisa Cl. Sci. (5)}, 8(3):469--507, 2009. 

\bibitem{Scheven-Schmidt:2010}
C.~Scheven and T. Schmidt,
\newblock Asymptotically regular problems. I. Higher integrability. 
\newblock {\em J. Differential Equations}, 248(4):745--791, 2010.

\end{thebibliography}

\end{document}